\documentclass[10pt]{article}

\usepackage{enumerate}
\usepackage{graphicx}
\usepackage{amsmath}
\usepackage{amssymb}
\usepackage[latin1]{inputenc}
\def\r{\mathbb{R}}
\def\b{\mathbb{B}}
\def\d{\mathbb{D}}
\def\n{\mathbb{N}}
\def\c{\mathbb{C}}
\def\s{\mathbb{S}}

\def\nor{{\cal N}}
\def\pro{{\cal P}}
\newcommand{\df}{ \stackrel{\rm def}{=}}
\newcommand{\cte}{\small \text{const.}}

\newcommand{\metri}[1]{{\cal S}_{#1}}

\usepackage{amsthm}

\newtheorem{teorema}{Theorem}
\newtheorem*{teoremaintro}{Theorem}
\newtheorem*{teoremaintroA}{Theorem A}
\newtheorem*{teoremaintroB}{Theorem B}

\newtheorem{claim}{Claim}
\newtheorem{lema}{Lemma}
\newtheorem{corolario}{Corollary}
\newtheorem{proposicion}{Proposition}
\newtheorem{remark}{Remark}

\newcommand{\dist}{\operatorname{dist}}
\newcommand{\re}{\operatorname{Re}}
\newcommand{\im}{\operatorname{Im}}

\newcommand{\met}[1]{\operatorname{\cal S}_{#1}}
\newcommand{\eucli}{\left<\cdot,\cdot\right>}
\newcommand{\intc}{\operatorname{Int}}

\newcommand{\longui}{\operatorname{len}}

\newfont{\goti}{cmr10 at 11pt}
\newcommand{\itx}[2]{{\goti(#1}$_{#2}${\goti)}}

\begin{document}

\title{Complete proper minimal surfaces in convex bodies of $\r^3$}
\author{\\ Francisco Martín\thanks{Research partially supported by MCYT-FEDER Grant no. BFM2001-3489 \newline
} \and \\ Santiago Morales$^*$}
\date{\today}
\maketitle
\begin{abstract} 

Consider a convex domain $B$ of $\r^3$. We prove that there exist
complete minimal surfaces which are properly immersed in $B$. We also
demonstrate that if $D$ and $D'$ are convex domains with $D$ bounded
and the closure of $\overline{D}$ contained in $D'$ then  any minimal
disk whose boundary lies in the boundary of $D$, can be approximated
in any compact subdomain of $D$ by a complete minimal disk which is
proper in $D'$. We apply these results to study the so called type
problem for a minimal surface: we demonstrate that the interior of any
convex region of $\r^3$ is not a universal region for minimal surfaces, in the sense explained by Meeks and Pérez in \cite{mp}. 
\end{abstract}
\section{Introduction} \label{sec:intro}
The global theory of complete minimal surfaces in $\mathbb{R}^3$ has been
developed for almost two and a half centuries. One of the central
problems in this theory has been the Calabi-Yau problem, which we now
introduce.  By the maximum principle for harmonic functions, there are
no  compact minimal surfaces in $\mathbb{R}^3$. Moreover, a basic observation of the classical examples of complete nonflat minimal surfaces (catenoid, helicoid, Riemann minimal examples, ...) reveals that one cannot bound any coordinate function for these surfaces. 
Even more, none of these examples is contained in a halfspace. This facts motivated E. Calabi to conjecture that there were no complete bounded minimal surfaces in $\r^3$, and later on S. T. Yau sharpened  Calabi's conjecture to ask if one can find a complete nonflat minimal surface in a halfspace of $\r^3$. This Calabi-Yau question is an extremely active field of research nowadays: very recently, T. H. Colding and W. P . Minicozzi [1] have answered the question in the negative when the surface is assumed to be embedded and with finite topology (other partial answers to the Calabi-Yau's question in the embedded setting have been given by Meeks and Rosenberg \cite{helicoid} and by Meeks, Pérez and Ros \cite{mprI}).

Embeddedness creates a dichotomy in the Calabi-Yau's question. The first example of a complete minimal surface with a bounded coordinate function was a disk constructed by L. P. Jorge and F. Xavier in 1981 [5]. The analytic arguments introduced by Jorge and Xavier in their construction were quite ingenious, and have been present in almost all the papers devoted to find complete minimal surfaces with some kind of boundedness on their coordinates. For instance, their idea of using a labyrinth of compact sets around the boundary of a simply connected domain, jointly with Runge's theorem in order to get completeness, was afterwards used by N. Nadirashvili in a more elaborated way to construct an example of a complete minimal surface in a ball of $\r^3$.
\begin{teoremaintro}[Nadirashvili, \cite{nadi}]
There exists a complete minimal immersion $f \colon \mathbb{D}
\rightarrow \b$ from the open unit disk $\d$ into the open unit ball $\b \subset \r^3$. Furthermore the immersion can be constructed with negative Gaussian curvature.
\end{teoremaintro}

However, Nadirashvili's technique did not guarantee the immersion $f
\colon \mathbb{D} \rightarrow \b$ was 
proper, where by proper we mean that $f^{-1}(C)$ is compact for any $C
\subset \b$ compact. Recently, Martín and Morales \cite{tran} introduced an
additional ingredient into Nadirashvili's machinery in order to
produce a complete minimal disk which is properly immersed in a ball
of $\r^3$. Similar ideas had been also used by Morales \cite{propi} to
construct a minimal disk which is properly immersed in $\mathbb{R}^3$. An
example of a complete proper minimal annulus which lies between two
parallel planes was constructed earlier by H. Rosenberg and E. Toubiana in \cite{rotu}. This example is related to the previous construction of Jorge and Xavier.

In this paper we answer the question to the existence of complete
simply connected minimal surfaces which are proper in convex regions of space. To be more precise, our main Theorem (see Theorem \ref{th:refinitivo} and Corollary \ref{co:esperanza}) establishes that:
\begin{teoremaintroA}
If $ B \subset \r^3$ is a convex domain (not necessarily bounded or smooth), then there exists a complete proper minimal immersion $\psi:\d \rightarrow B$.
\end{teoremaintroA}

The proof of this theorem is divided in a two-stage process. First we
prove the theorem for bounded convex regular domains, i.e. domains
whose boundary is a compact analytic surface of $\r^3$ (Theorem
\ref{th:mari}.) Later, we make use of a classical result by Minkowski
(Theorem \ref{th:minko}) to approximate in terms of Hausdorff
distance the interior of an arbitrary convex region $B$ by an increasing
sequence of bounded convex regular domains $\{B^n\}.$
Theorem~\ref{th:mari}  give us the existence of a complete minimal
disk $M_n$ which is  properly immersed in $B^n$. Roughly speaking, the
desired minimal disk $M$ is the limit of  the family
$\{M_n\}$. Concerning this point, we are indebted to W. H. Meeks who
suggested to us the way of making use of  Theorem \ref{th:mari} to
construct the above sequence of minimal disks to get a good limit
$M$. Item (b) in Theorem \ref{th:mari} below restricts the behavior of the above
approximation  in a way that can  be used in the proof of the properness for the limit immersion. 
It is important to emphasize that we  not only establish the existence
of properly immersed  minimal disks in a convex domain, but we also prove that: 
\begin{teoremaintroB}
Let $D'$ be a convex domain of space, and consider $D$ a bounded
convex regular domain, with $\overline D \subset D'$. Then, any
minimal disk whose boundary lies in $\partial D$, can be approximated
in any compact subdomain of $D$ by a complete minimal disk which is {\em proper} in $D'$. 
\end{teoremaintroB}    

The lemmas we will use in proving Theorem \ref{th:mari} are stated and
proved in Section \ref{sec:lemmata}. Lemma \ref{lem:propia} is used to
get Lemma \ref{lem:nadi} and it represents the main ingredient used to obtain properness. It essentially asserts that a minimal disk with boundary can be perturbed  outside a compact set in such a way that the boundary of the resulting surface achieves the boundary of a prescribed convex domain. 
 
 Lemma \ref{lem:nadi} is the main tool in the proof of Theorem \ref{th:mari}. It will allow us to construct, in a recursive way, a  convergent sequence of compact minimal disks used to prove this theorem. It can be considered to be a refinement of the corresponding Nadirashvili's lemma in \cite{nadi}. However, we must adapt his ideas about balls to the setting of general convex bodies of space. Hence, our deformation increases the intrinsic diameter, but it controls the extrinsic geometry in $\r^3$, in such a way that the perturbed surface is included in a prescribed convex set. Moreover, in order to obtain the boundary behavior described in items (b.3) and (b.4) of Lemma \ref{lem:nadi}, we have to apply Lemma \ref{lem:propia} in the way we have explained in the previous paragraph. These properties will be crucial in the delicate argument that gives properness for the limit immersion. 

In some sense, our Theorem A is related with an intrinsic question
associated to the underlying complex structure: the so called {\em
  type problem} for a minimal surface $M$, i.e. whether $M$ is
hyperbolic or parabolic (as we have already noticed, the elliptic
(compact) case is not possible for a minimal surface). Classically, a
Riemann surface without boundary is called {\em hyperbolic} if it
carries a nonconstant positive superharmonic function, and {\em
  parabolic} if it is neither compact nor hyperbolic. In the case of a
Riemann surface with boundary, we say that $M$ is {\em parabolic} if every
bounded harmonic function on $M$ is determined by its boundary values,
otherwise $M$ is called {\em hyperbolic.} It turns out that the parabolicity
for Riemann surfaces without boundary is equivalent to the recurrence
of Brownian motion on such surfaces. If the boundary of $M$ is
nonempty, then $M$ is called {\em parabolic} if, and only if, there exists a point
$p$ in the interior of $M$ such that the probability of a Brownian
path beginning at $p$, of hitting the boundary $\partial M$ is $1$ (see \cite{grigo} for more details.) 

In this setting, given a connected region $W \subset \r^3$  which is either open or the closure of an open set, we say that $W$ is {\em universal for surfaces} if every complete, connected, properly immersed minimal surface $M \subset W$ is either recurrent ($\partial M=\emptyset$) or a parabolic surface with boundary. The open question of determining which  regions of space are universal for surfaces has been proposed by W. H. Meeks and  J. Pérez in \cite{mp}. 
Theorem A implies that a convex domain of $\r^3$  is not {\em universal for surfaces.} In contrast with this result, it is known \cite{ckmr} that the closure of a convex domain is universal for surfaces.

\section{Preliminaries and Notation} \label{sec:pre}

Our objective in this section is to briefly summarize the notation and
results about minimal surfaces and convex geometry that we will use in the paper.

\subsection{ Minimal surface background}
The theory of complete minimal surfaces is closely related to complex
analysis of one variable. This is due to the fact that any such surface is given by a triple $\Phi=(\Phi_1, \Phi_2, \Phi_3)$ of holomorphic 1-forms defined on some Riemann surface such that:
\begin{equation} \label{eq:conforme}
\Phi_1^2+\Phi_2^2+\Phi_3^2=0;
\end{equation}
\begin{equation} \label{eq:bilbao}
\|\Phi_1\|^2+\|\Phi_2\|^2+\|\Phi_3\|^2 \neq 0;
\end{equation}
and all periods of the $\Phi_j$ are purely imaginary; here we consider
$\Phi_i$ to be a holomorphic function times $dz$ in a local parameter $z$. Then the minimal immersion $X:M \rightarrow \r^3$ can be parametrized by $z \mapsto \mbox{Re} \int^z \Phi.$
The above triple is called the Weierstrass representation of the immersion $X$. Usually, the first requirement (\ref{eq:conforme}) (which ensures the conformality of $X$) is guaranteed by introducing the formulas:
$$\Phi_1 =\tfrac12 \left( 1-g^2\right) \, \eta, \quad \Phi_2 =\tfrac{\rm i}2 \left( 1+g^2\right) \, \eta, \quad \Phi_3= g \, \eta, $$
with a meromorphic function $g$ (the stereographic projection of the Gauss map) and holomorphic 1-form $\eta$. In this article, all the minimal immersions are defined on simply connected domains of $\c$. Then, the Weierstrass 1-forms have no periods, and so the only requirements are (\ref{eq:conforme}) and (\ref{eq:bilbao}). In this case, the differential $\eta$ can be written as $\eta=f(z) \, dz$. 
The metric of $X$ can be expressed as
\begin{equation} \label{eq:metric}
{\metri X}^2=\tfrac12 \|\Phi\|^2=\left(\tfrac12\left(1+|g|^2\right) |f| |dz|\right)^2.
\end{equation}
Throughout the paper, we will use several orthonormal bases of $\r^3$. Given $X:\Omega \rightarrow \r^3$ a minimal immersion and $S$ an orthonormal basis, we will write the Weierstrass data of $X$ in the basis $S$ as $${\Phi}_{(X,S)}=(\Phi_{(1,S)},\Phi_{(2,S)},\Phi_{(3,S)}), \quad f_{(X,S)}, \quad g_{(X,S)}, \quad \eta_{(X,S)}.$$
Similarly, given $v \in \r^3$, we will let $v_{(k,S) }$ denote the $k$-th coordinate of $v$ in $S$. The first two coordinates of $v$ in this basis will be represented by $v_{(*,S)}=\left(v_{(1,S)},v_{(2, S)} \right)$.

Given a curve $\alpha$ in $\Omega$, by $\longui(\alpha, \metri X)$ we
mean the length of  $\alpha$ with respect to the metric $\metri X$. Given a subset $W \subset \Omega$, we define:
\begin{itemize}
\item $\dist_{(W,\metri X)}(p,q)=\inf \{\longui(\alpha, \metri X) \: | \: \alpha:[0,1]\rightarrow W, \; \alpha(0)=p,\alpha(1)=q \}$, for any $p,q\in W$;
\item $\dist_{(W,\metri X)}(T_1,T_2)=\inf \{\dist_{(W,\metri X)}(p,q) \;|\;p \in T_1, \;q \in T_2 \}$, for any $T_1, T_2 \subset  W$;
\item $\text{diam}_{\metri X}(W)=\sup \{\dist_{(W,\metri X)}(p,q) \;|\;p,q \in W \}$. 
\end{itemize}
The Euclidean metric on $\c$ will be denoted as $\langle \cdot,\cdot \rangle$. Note that ${\metri X}^2=\lambda_X^2 \,\langle \cdot,\cdot \rangle$, where the conformal coefficient $\lambda_X$ is given by (\ref{eq:metric}). 
 
Given a domain $D\subset \c$, we will say that a function, or a 1-form, is harmonic, holomorphic, meromophic, ... on $\overline D$, if it is harmonic, holomorphic, meromorphic, ... on a domain containing $\overline D$.

Let $P$ be a simple closed polygonal curve in $\c$. We let $\intc P$
denote the bounded connected component of $\c \setminus P.$ We will assume that the origin is in the interior region of all the polygons that appears in the paper.  Given
$\xi>0$, small enough, we define $P^\xi$ to be the parallel polygonal
curve in $\intc P$, satisfying the property that the distance between parallel sides is equal to $\xi$. Whenever we write $P^\xi$ in the paper we are assuming that $\xi$ is small enough to  define the polygon properly.

\subsubsection{The López-Ros transformation}
The proof of Lemmas \ref{lem:propia} and \ref{lem:nadi} exploits what has come to be called the López-Ros transformation. If $(g,f)$ are the Weierstrass data of a minimal immersion $X:\Omega \rightarrow \r^3$ ($\Omega \subset \c$ simply connected), we define on $\Omega$ the data
\begin{equation} \label{eq:mambru}
\widetilde g= \frac{g}{h}, \qquad \widetilde f= f \cdot h,
\end{equation}
where $h:\Omega \rightarrow \c$ is a holomorphic function without
zeros. Notice that the new meromorphic data satisfies
(\ref{eq:conforme}) and (\ref{eq:bilbao}). So, the new data defines a minimal immersion $\widetilde X: \Omega \rightarrow \r^3$. This method provides us with a powerful and natural tool for deforming minimal surfaces. From our point of view, the most important property of the resulting surface is that the {\em third coordinate function is preserved.}  Note that the intrinsic metric is given by (\ref{eq:metric}) as
 \begin{equation} \label{eq:metric2}
{\metri{\widetilde X}}^2=\left(\tfrac12\left(|h|^2+|g|^2\right)\, |f|\,|dz| \right)^2.
\end{equation}
This means that we can increase the intrinsic distance in a prescribed compact $K \subset \Omega$, by using  suitable functions $h$. These functions will be provided by Runge's theorem.
\subsection{Background on convex geometry}
Convex geometry is a classical subject with a large literature. To make this article self-contained, we will describe the concepts and results we will need. 
A convex set of $\r^n$ with nonempty interior is called {\em a convex body}. A theorem of H. Minkowski (cf. \cite{mingorebulgo}) states that every convex body $C$ in $\r^n$ can be approximated (in terms of Hausdorff metric) by a sequence $C_k$ of `analytic' convex bodies.
\begin{teorema} \label{th:minko}
Let $C$ be a convex body in $\r^n$. Then there exists a sequence $\{C_k \}$ of convex bodies with the following properties
\begin{enumerate}[\rm 1.]
\item $C_k \searrow C$;
\item $\partial C_k$ is an analytic $(n-1)$-dimensional manifold;
\end{enumerate}
\end{teorema}
A modern proof of this result can be found in \cite[\S 3]{meeksyau}. 

Given $E$ a bounded regular convex domain of $\r^3$ and $p \in
\partial E$, we will let $ \kappa_2(p) \geq 0$ denote the largest
principal curvature of $\partial E$ at $p$ (associated to the inward
pointing unit normal.) Moreover, we write $$\kappa_2(\partial E) \df \mbox{max} \{ \kappa_2(p) \: : \: p \in \partial E \}.$$
If we consider $\nor: \partial E \rightarrow \s^2$ the outward
pointing unit normal or Gauss map of $\partial E$, then there exists a constant $a>0$ (depending on $E$) such that $\partial E_t=\{p+ t\cdot \nor(p) \; : \; p \in \partial E\}$ is a regular (convex) surface $\forall t \in [-a, +\infty[$. We label $E_t$ as the convex domain bounded by $\partial E_t$. The normal projection to $E$ is represented as 
$$\pro_E:\r^3\setminus E_{-a}\longrightarrow \partial E,$$
$$p+t \cdot \nor(p) \mapsto p.$$
Finally, we define the `extended' Gauss map $\nor_E : \r^3 \setminus E_{-a} \longrightarrow \s^2$ as $\nor_E(x)=\nor(\pro_E(x)).$

\section{Preliminary Lemmas} \label{sec:lemmata}
As we have indicated in the introduction, the proofs of the main
results of the paper require the following two technical lemmas.
Actually, Lemma \ref{lem:propia} is used in proving Lemma \ref{lem:nadi}.
\begin{lema} \label{lem:propia}
Let $E$ and $E'$ be two  regular bounded convex domains in $\r^3$, with $0\in E \subset \overline{E} \subset E'$.  Let $X:O\longrightarrow\r^3$ be a conformal minimal immersion defined on a simply connected domain $O$, $0\in O$, with $X(0)=0$. Consider a polygon $ P$ with $ P\subset  O$,  satisfying:
\begin{equation}\label{eslabon}
 X(O\setminus \intc P)\subset E'\setminus \overline{E}.
\end{equation}
Then, for any $ b_1, b_2>0$, such that $E'_{-b_2}$ and $E_{-2 b_2}$ exist, there exist a polygon $Q$ and a conformal minimal immersion $Y:\overline{\intc Q}\longrightarrow \r^3$, with $Y(0)=0$, such that:
\begin{enumerate}[\rm ({a.}1)]
\item $ P\subset \intc Q\subset\overline{\intc Q}\subset  O$;
\item $\|Y(z)- X(z)\|< b_1$, $\forall z\in\overline{\intc  P}$;
\item $Y(Q)\subset E'\setminus E'_{- b_2}$;
\item $Y(\intc Q\setminus\intc  P)\subset \r^3\setminus E_{-2 \, b_2}$.
\end{enumerate}
\end{lema}

\begin{lema} \label{lem:nadi}
Let $E$ and $E'$ be two  regular bounded convex domains in $\r^3$, with $0\in E \subset \overline{E} \subset E'$. Let $P$ be a polygon, $X: \overline{\intc P} \longrightarrow\r^3$ be a conformal minimal immersion, with $X(0)=0$, and $\varepsilon$, $a$ and $b$ positive constants, such that:
\begin{enumerate}[\rm 1)]
\item $X(\overline{\intc P\setminus\intc P^\varepsilon} )\subset E\setminus E_{-a}$;
\item $\sqrt{\left(a+\frac1{\kappa_2(\partial E)}\right)^2+a^2}-\frac1{\kappa_2(\partial E)}+\varepsilon<\frac{\dist(\partial E, \partial E')}2.$ 
\end{enumerate}
Then, there exist a polygon $Q$ and a conformal minimal immersion $Y: \overline{\intc Q} \rightarrow \r^3$, with $Y(0)=0$, and verifying:
\begin{enumerate}[\rm ({b}.1)]
\item $\overline{\intc P^\varepsilon} \subset \intc Q \subset \overline{\intc Q} \subset \intc P$;
\item $\sigma<\dist_{(\overline{\intc Q},\met Y)}(z,P^\varepsilon)$, $\forall z\in Q$, where $\sigma>0$ is given by the equation \\$\sqrt{\left(a+\frac1{\kappa_2(\partial E)}\right)^2+\left(2 \sigma+a\right)^2}-\frac1{\kappa_2(\partial E)}+\varepsilon=\frac{\dist(\partial E, \partial E')}2$;
\item $Y(Q) \subset E' \setminus E'_{-b}$;
\item $Y(\intc Q \setminus \intc P^\varepsilon) \subset \r^3 \setminus E_{-2 (a+b)}$;
\item $\| Y(z)-X(z)\| <\varepsilon$, $\forall z\in \overline{\intc P^\varepsilon}$.

\end{enumerate} 
\end{lema}

\subsection{Proof of Lemma \ref{lem:propia}}
The proof of Lemma \ref{lem:propia} consists of deforming the original immersion $X$ outside $\overline{\intc P}$. The idea is to make this a two-stage process. 
In the first stage, we make use of meromorphic functions with single poles to apply successive López-Ros transformations. These poles $p_1, \ldots , p_n$ are distributed around the polygon $P$. In this way, we deform our original surface in a neighborhood of the points $p_i$, following  the direction of $\nor_E(X(p_i))$, $i=1, \ldots , n$. In a second stage,
we consider suitable curves joining two consecutive poles, $p_i$ and $p_{i+1}$, and we deform the surface again  along these curves. This time, the deformation acts tangentially to $\partial E$, by using Runge's functions as López-Ros parameters. These functions achieve very big values on  the above mentioned curves. In both stages, the functions that appear in López-Ros transformations are close to $1$ in a neighborhood of $\overline{\intc P}$. It guarantees that the resulting immersion $Y$ will be close to $X$ in $\overline{\intc P}$.

Along this proof, we denote $\nor$ and $\pro$ as the extended Gauss map associated to $E$ and the normal projection to $\partial E$, respectively. We also define  two important constants that are chosen as follows:
\begin{itemize}\label{defconstant}
\item $\mu=\max \{ \dist(p, \partial E) \; : \; p \in E'\}$;
\item $\epsilon_0>0$ is taken small enough to satisfy all the inequalities appearing in this section. This choice will only depend on the data of the lemma.
\end{itemize}
\subsubsection{The first deformation stage}
We start by  choosing the aforementioned points $p_i$, $i=1, \ldots, n$. 
\begin{claim} \label{frame}
There exist a simply connected domain $W$, with  $\overline{\intc P}\subset W\subset \overline W\subset O$, and a set of points $\{ p_1,p_2, \ldots,p_n \}$ included in $W\setminus \overline{\intc P}$, satisfying the following properties:
\begin{enumerate}[1.]
\item The segments $\overline{p_1\,p_2},\ldots,\overline{p_{n-1}\, p_n},\overline{p_n\,p_{n+1}}$ form a polygon $\widehat P$ in $W\setminus\overline{\intc P}$ (we adopt the convention $p_{n+1}=p_1$);
\item There exist open disks $B^i\subset W\setminus\overline{\intc P}$ satisfying  $p_i,p_{i+1}\in B^i$, for all  $i\in\{1,\ldots,n\}$, and 
\begin{equation}\label{tarariX}
\|X(z)-X(w)\|<\epsilon_0, \quad \forall z, w \in B^i;
\end{equation}
\item For each  $p_i$, $i=1, \ldots , n$, there exists an orthonormal basis of $\r^3$, $S_i=\{e_1^i,e_2^i,e_3^i\}$, with $e_1^i=\nor(X(p_i))$, and satisfying:
\begin{equation}\label{tarariS}
\left\|e_j^i-e_j^{i+1}\right\|<\frac{\epsilon_0}{3\mu},\quad \forall j \in \{1,2,3\},
\end{equation}
and
\begin{equation}\label{tararif}
f_{(X,S_i)}(p_i)\not=0;
\end{equation}
\item For  each $p_i$, $i=1, \ldots n$,  there exists a complex number $\theta_i$, satisfying $|\theta_i|=1$, $\im \theta_i\not=0$, and:
\begin{equation} \label{borringa}
\left|\frac{\theta_if_{(X,S_i)}(p_i)}{|f_{(X,S_i)}(p_i)|}- 1 \right| < \frac{\epsilon_0}{3\mu},
\end{equation}
\end{enumerate}
\end{claim}
\begin{remark}
It is important to note that Properties (\ref{tarariX}) and (\ref{tarariS}) are cyclic: i.e. they are true for $i=n$ considering $p_{n+1}=p_1$, $S_{n+1}=S_1$, and $B^{n+1}=B^1$.
\end{remark}

\begin{proof}
As $P$ is a piecewise regular curve, we know that $\nor(X(P))$ omits an open set $U_0$ of $\s^2$. Hence, we can find a simply connected domain $W$, verifying $\overline{\intc P} \subset W$ and ${\cal N}(X(W \setminus \intc P)) \subset \s^2 \setminus U_0$. Let $V_1$ and $V_2$ be a smooth orthonormal basis of tangent vector fields on $\s^2 \setminus U_0$. Then, we define:
$ \hat e_1(z)={\cal N}(X(z)),$ $ \hat e_2(z)=V_1\left({\cal N}(X(z))\right),$ $\hat e_3(z)=V_2\left({\cal N}(X(z))\right)$, $z\in W\setminus\overline{\intc P}$.
\begin{figure}[h]
\begin{center}
\includegraphics[width=.55 \textwidth]{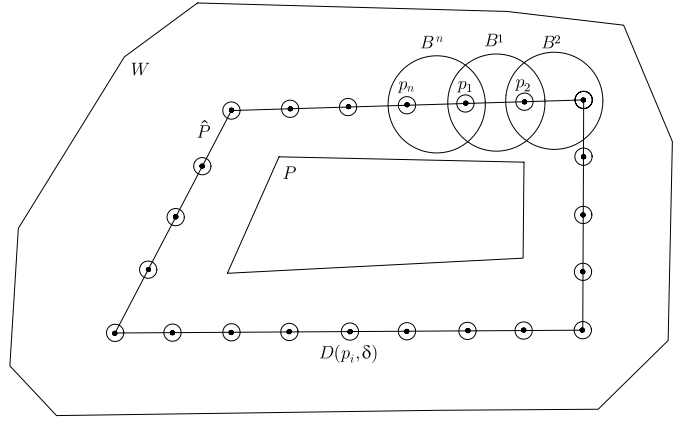}
\caption{The points $p_i$.}\label{lospe}
\end{center}
\end{figure}

The points $p_1,\ldots,p_n$ will be taken in $W\setminus\overline{\intc P}$ to satisfy Statements {\em 1}, {\em 2}, and the following property:
\begin{equation}\label{tarariquetevi}
\left\|\hat e_j(p_i)-\hat e_j(p_{i+1})\right\|<\epsilon_0/2,\quad \forall j \in \{1, 2,3\}, \quad \forall i \in \{1, \ldots, n\}.
\end{equation}
If the points $p_i$, $i=1, \ldots, n$, are close enough (i.e.,  the number of points is large enough), then 
these properties are a direct consequence of the uniform continuity of $X$ and the fields $\hat e_j$, respectively. 
In order to obtain Statement {\em 3}, we consider the Gauss map $G$ of $X$. We can write $G(p_i)= \sum_{j=1}^3 \hat g_j^i \cdot \hat e_j(p_i), \quad \hat g_j^i \in [0,1].$
Take $a \in [0,1] \setminus \{ \hat g_2^i, i=1, \ldots , n\}$, and define $e_1^i=\hat e_1(p_i)$, $e_2^i=-\sqrt{1-a^2} \hat e_2(p_i)+a \hat e_3(p_i)$ and $e_3^i=a \hat e_2(p_i)+\sqrt{1-a^2}\hat e_3(p_i)$. Then, (\ref{tarariS}) is a direct consequence of (\ref{tarariquetevi}). Moreover, note that $e_3^i \neq G(p_i)$, $\forall i$, and so (\ref{tararif}) trivially holds. 
Finally, the choice of the complex numbers $\theta_i$, $i=1, \ldots , n$, is straightforward. 
\end{proof}

\begin{figure}[hbtp]
\begin{center}
\includegraphics[width=4cm]{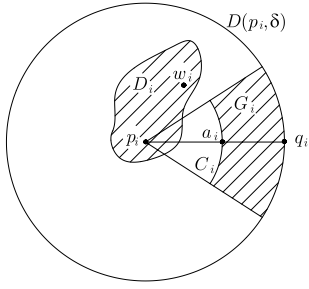}
\caption{The disk $D(p_i,\delta)$.}\label{que}
\end{center}
\end{figure}

A  small enough real number $\delta$ is chosen in $]0,1[$ and verifying the following list of properties, (see Figure \ref{lospe}):
\begin{enumerate}[ ({A}1)]
\item $\overline{\intc \widehat P\setminus \cup_{k=1}^n D(p_k,\delta)}$ is simply connected, where $D(p_k,\delta)$ means the disk centered at $p_k$ with radius $\delta$;
\item $\overline{D(p_i,\delta)\cup D(p_{i+1},\delta)}\subset B^i,\quad \forall i=1,\ldots,n$;
\item $D(p_i,\delta)\cap D(p_k,\delta)=\emptyset,\quad \forall i=1,\ldots,n$, and $k\not=i$;
\item $\delta\max_{\overline{D(p_i,\delta)}}\{|f_{(X,S_i)}|\} <2\epsilon_0,\quad \forall i=1,\ldots,n$;
\item $\delta \cdot \frac{\max_{\overline{D(p_i,\delta)}} \{|f_{(X,S_i)}g^2_{(X,S_i)}|\}}{|\im\theta_i|} <2\epsilon_0,\quad \forall i=1,\ldots,n$;
\item $\delta \cdot \max_{\overline{D(p_i,\delta)}} \{\|\phi^0\|\}<\epsilon_0,\quad \forall i=1,\ldots, n$;
\item $\displaystyle 3 \mu \, \frac{\max_{w \in \overline{D(p_i,\delta)}} \{|f_{(\Phi^0,S_i)}(w)-f_{(\Phi^0,S_i)}(p_i)|}{|f_{(\Phi^0,S_i)}(p_i)|}<\epsilon_0,\quad \forall i=1,\ldots, n.$
\end{enumerate}
If we label ${\cal E}=\overline{\intc \widehat P\setminus \cup_{k=1}^n D(p_k,\delta)}$, then we define the constant $l$ as:
\begin{equation}\label{quefrio}
l \df \sup_{z\in{\cal E}} \left\{\dist_{\left({\cal E},\eucli\right)}(0,z) \right\}+2\pi\delta+\delta+1.
\end{equation}
We are now ready to construct a sequence of Weierstrass representations that will provide us the immersion we are looking for in this deformation stage. To be more precise, we will construct a sequence $\Psi_1,\ldots,\Psi_n$ in a recursive process, where the element $\Psi_i=\{\Phi^i,k_i,a_i,C_i,G_i,D_i\}$ is composed of:
\begin{itemize}\label{alcatel}
\item $k_i$ a positive constant;
\item $\Phi^i$ is a Weierstrass representation defined on $\overline W$. The points $p_1,\ldots,p_i$ are poles of  $\Phi^i$; and the points $w_j=p_j-k_j\theta_j$, $j=1, \ldots i$, are either zeros or poles of  $\Phi^i$. We write $\Phi^i=\phi^idz$ where  $\phi^i:\overline W\rightarrow \c^3$ is a meromorphic map;
\item $a_i$ is a point in the segment $\overline{p_iq_i}$, where $q_i=p_i+\delta\in\partial D(p_i,\delta)$;
\item $C_i$ is an arc of a circumference centered at  $p_i$ and containing $a_i$;
\item $G_i$ is  a closed annular sector bounded by $C_i$, a piece of $\partial D(p_i,\delta)$ and two radius of these circumferences, as Figure \ref{que} indicates;
\item $D_i$ a simply connected domain of $\c$ verifying $\overline{D_i}\cap G_i=\emptyset$ and
$\{p_i,w_i\}\subset D_i\subset\overline{D_i}\subset D(p_i,\delta).$
\end{itemize}

\begin{remark}\label{n+1}
In what follows, we will adopt the convention that $\Psi_{n+1}= \Psi_1$. Therefore, we will write  $\Phi^{n+1}=\Phi^1$, $k_{n+1}=k_1$, $\ldots$, $p_{n+1}=p_1$, $q_{n+1}=q_1$, $\ldots$
\end{remark}
The sequence $\Psi_1,\ldots,\Psi_n$ is constructed to satisfy:
\begin{enumerate}[ ({B}1$_{i}$)]
\item $\delta\max_{\overline{D(p_k,\delta)}}\{|f_{(\Phi^i,S_k)}|\}<2\epsilon_0,\quad k=i+1,\ldots,n$;
\item $\displaystyle \delta\frac{\max_{\overline{D(p_k,\delta)}} \{|f_{(\Phi^i,S_k)}g^2_{(\Phi^i,S_k)}|\}}{|\im\theta_k|}<2\epsilon_0,\quad k=i+1,\ldots,n$;
\item $\displaystyle 3 \mu \, \frac{\max_{w \in \overline{D(p_k,\delta)}} \{|f_{(\Phi^i,S_k)}(w)-f_{(\Phi^0,S_k)}(p_k)|}{|f_{(\Phi^0,S_k)}(p_k)|\}}<\epsilon_0,\quad k=i+1,\ldots,n$;
\item $\|\re\int_{\alpha_z}\Phi^i\|<\epsilon_0,\quad \forall z\in C_i$, where $\alpha_z$ is a piece of $C_i$ connecting $a_i$ with $z$;
\item $\Phi_{(3,S_i)}^i=\Phi_{(3,S_i)}^{i-1}$, where $\Phi_{(k,S_i)}^i$ represents the $k$-th coordinate of the triple $\Phi^i$ in the frame $S_i$;
\item $\|\phi^i(z)-\phi^{i-1}(z)\|<\frac{\epsilon_0}{4n l},\quad \forall z\in\overline W\setminus (D(p_i,\delta)\cup(\cup_{k=1}^{i-1} D_k))$;
\item $\|\re\int_{\overline{q_ia_i}}\Phi^i-\re \int_{\overline{q_{i-1}a_{i-1}}}\Phi^{i-1}\|< 13 \epsilon_0$, (for $i=2,\ldots,n+1$).
\item For all $z\in G_i$ one has
$$\left\|\left(\re\int_{\overline{q_iz}}\Phi_1^i\right)e_1^i+\left(\re\int_{\overline{q_iz}}\Phi_2^i\right)e_2^i\right.
-\left.\tfrac12\left(\int_{\overline{q_iz}}\frac{k_idw}{w-p_i}\right)|f_{(\Phi^0,S_i)}(p_i)|\:e_1^i \right\|<4\epsilon_0.$$

\end{enumerate}
The above properties are true for $i=1,\ldots,n$, except for \itx{B1}{i}, \itx{B2}{i} and \itx{B3}{i} that only hold for $i=1,\ldots,n-1$. In the same way, Property \itx{B7}{i} is only valid for $i=2,\ldots,n+1$ (see Remark \ref{n+1}).
Observe that Properties \itx{B5}{i} and \itx{B8}{i} tell us that the deformation of our surface around the points $p_i$ follows the direction of $e_1^i=\nor(X(p_i))$. 

We define $\Psi_1,\ldots,\Psi_n$ in a recursive way. Let $\Phi^0=\phi^0dz$ be the Weierstrass representation of the immersion $X_0=X$. We denote $\Psi_0=\{\Phi^0\}$. Assume we have defined $\Psi_{i-1}$ verifying Properties \itx{B1}{i-1}$,\ldots,$\itx{B8}{i-1}. To construct $\Psi_i$ starting from $\Psi_{i-1}$ we will use Properties
\itx{B1}{i-1}, \itx{B2}{i-1} and \itx{B3}{i-1}. In the case $i-1=0$, these properties are a consequence of (A4), (A5) and (A7).

The Weierstrass data $\Phi^i$, in the orthogonal frame $S_i$, are determined by the López-Ros transformation:
$$f_{(\Phi_i,S_i)}=f_{(\Phi_{i-1},S_i)} \cdot h_i, \quad g_{(\Phi_i,S_i)}=\frac{g_{(\Phi_{i-1},S_i)}}{h_i},$$
where $h_i(z)=\frac{k_i\theta_i}{z-p_i}+1$. The constant $k_i>0$ is taken small enough to satisfy Properties \itx{B1}{i}, \itx{B2}{i}, \itx{B3}{i} and \itx{B6}{i}. Notice that this is possible thanks to the fact $\Phi^i \stackrel{k_i\rightarrow 0}{\longrightarrow} \Phi^{i-1}$, uniformly on $\overline W\setminus (D(p_i,\delta)\cup(\cup_{k=1}^{i-1} D_k))$. 
Furthermore, Property \itx{B5}{i} trivially follows from the definition of $\Phi^i$.

The point $a_i$ is the first point in the (oriented) segment $\overline{q_ip_i}$, such that:
\begin{equation}\label{IA}
\tfrac12|f_{(\Phi^0,S_i)}(p_i)|\int_{\overline{q_ia_i}} \frac{k_idw}{w-p_i}=3 \mu.
\end{equation}
We take $D_i$ to be a simply connected domain containing the pole and the zero, $w_i=p_i-k_i\theta_i$, of $h_i$ and verifying  $\overline{D_i}\subset D(p_i,\delta)$ and $\overline{D_i}\cap\overline{q_ia_i}=\emptyset$. It is possible because $\im\theta_i\not=0$, and so $w_i\not\in \overline{p_iq_i}$.

The next step is to prove Property \itx{B8}{i}. We will use complex notation, that is to say $a+{\rm i}\, b\equiv ae_1^i+be_2^i$. Consider $z\in\overline{q_ia_i}$. Taking (\ref{borringa}) and (\ref{IA}) into account one obtains:
\begin{multline*}
\left|\left(\re\int_{\overline{q_iz}}\Phi_{(1,S_i)}^i\right)+{\rm i} \left(\re\int_{\overline{q_iz}}\Phi_{(2,S_i)}^i\right)-
\tfrac12\left(\int_{\overline{q_iz}}\frac{k_idw}{w-p_i}\right)|f_{(\Phi^0,S_i)}(p_i)| \right|<\\
\left| \left(\re\int_{\overline{q_iz}}\Phi_{(1,S_i)}^i\right)+{\rm i}\left(\re\int_{\overline{q_iz}}\Phi_{(2,S_i)}^i\right)-  \tfrac12\left(\int_{\overline{q_iz}}\frac{k_idw}{w-p_i}\right)\overline{\theta_i f_{(\Phi^0,S_i)}(p_i)} \right|+ \epsilon_0=\\
\end{multline*}
From the expression of $\Phi^i$ and $h_i$, we have
\begin{multline*}\nonumber
=\tfrac12\left| \int_{\overline{q_iz}}\overline{f_{(\Phi^{i-1},S_i)}(w)\frac{k_i\theta_i}{w-p_i}dw}+\int_{\overline{q_iz}}\overline{f_{(\Phi^{i-1},S_i)}(w)dw}-\right.\\
\left.\int_{\overline{q_iz}}f_{(\Phi^{i-1},S_i)}(w)g^2_{(\Phi^{i-1},S_i)}(w)\frac{dw}{h_i(w)}-  \left(\int_{\overline{q_iz}}\frac{k_idw}{w-p_i}\right)\overline{\theta_i f_{(\Phi^0,S_i)}(p_i)} \right|+ \epsilon_0<\\
\tfrac12\left|\int_{\overline{q_iz}}\overline{\left(f_{(\Phi^{i-1},S_i)}(w)-f_{(\Phi^0,S_i)}(p_i)\right)\frac{k_i\theta_i}{w-p_i}dw} \;\right|+\tfrac12\left|\int_{\overline{q_iz}}\overline{f_{(\Phi^{i-1},S_i)}(w)dw}\right|+\\
\tfrac12\left|\int_{\overline{q_iz}}f_{(\Phi^{i-1},S_i)}(w)g^2_{(\Phi^{i-1},S_i)}(w)\frac{dw}{h_i(w)}\right|+ \epsilon_0<
\end{multline*}
by using (\ref{IA}), \itx{B1}{i-1}, \itx{B2}{i-1}, \itx{B3}{i-1} and the fact that $|h_i(w)|>|\im\theta_i|$, we can get upper bounds of the last addends,
$$<3\epsilon_0+\epsilon_0= 4\epsilon_0.$$
Thus, we have proved that Property \itx{B8}{i} holds for all $z\in \overline{p_i a_i}$. Therefore, if $C_i$ and $G_i$ are chosen sufficiently close to $a_i$ and $\overline{p_ia_i}$, respectively, we obtain Properties \itx{B4}{i} and \itx{B8}{i}.

To complete the construction of $\Psi_i$ we only need to check  \itx{B7}{i}. To do this, we write:
\begin{equation}\label{torres}
\begin{split}
\left\|\re\int_{\overline{q_ia_i}}\Phi^i-\re\int_{\overline{q_{i-1}a_{i-1}}}\Phi^{i-1}\right\|\leq & \left\|\sum_{k=1}^2\left(\re\int_{\overline{q_ia_i}}\Phi_{(k,S_i)}^i\right)e_k^i-\left(\re\int_{\overline{q_{i-1}a_{i-1}}}\Phi_{(k,S_{i-1})}^{i-1}\right)e_k^{i-1}\right\|+ \\
&\left\|\left(\re\int_{\overline{q_ia_i}}\Phi_{(3,S_i)}^i\right)e_3^i-\left(\re\int_{\overline{q_{i-1}a_{i-1}}}\Phi_{(3,S_{i-1})}^{i-1}\right)e_3^{i-1}\right\|.
\end{split}
\end{equation}

Using \itx{B8}{i}, (\ref{IA}), and (\ref{tarariS}), it is not hard to see:
\begin{equation}\nonumber
\left\|\sum_{k=1}^2\left(\re\int_{\overline{q_ia_i}}\Phi_{(k,S_i)}^i\right)e_k^i-\left(\re\int_{\overline{q_{i-1}a_{i-1}}}\Phi_{(k,S_{i-1})}^{i-1}\right)e_k^{i-1}\right\|\leq 9 \epsilon_0.
\end{equation}

Regarding the last addend of (\ref{torres}), we apply \itx{B5}{i-1} and \itx{B5}{i}, and obtain:
\begin{equation}\nonumber
\begin{split}
&\left\|\left(\re\int_{\overline{q_ia_i}}\Phi_{(3,S_i)}^i\right)\right.  \left.e_3^i- \left(\re\int_{\overline{q_{i-1}a_{i-1}}}\Phi_{(3,S_{i-1})}^{i-1}\right)e_3^{i-1}\right\|\leq\\
& \left|\re\int_{\overline{q_ia_i}}\Phi_{(3,S_i)}^{i-1}\right|+\left|\re\int_{\overline{q_{i-1}a_{i-1}}}\Phi_{(3,S_{i-1})}^{i-2}\right| \leq \delta\left(\max_{\overline{D(p_i,\delta)}}\{\|\phi^{i-1}\|\}+\max_{\overline{D(p_{i-1},\delta)}}\{\|\phi^{i-2}\|\}\right).
\end{split}
\end{equation}
Now, we use \itx{B6}{k}, for $k=1,\ldots,i-1$, and {(A6)} to deduce that:
\begin{equation}\nonumber
\left\|\left(\re\int_{\overline{q_ia_i}}\Phi_{(3,S_i)}^i\right)e_3^i- \left(\re\int_{\overline{q_{i-1}a_{i-1}}}\Phi_{(3,S_{i-1})}^{i-1}\right)e_3^{i-1}\right\|
\leq 2\epsilon_0+2\delta\epsilon_0<4\epsilon_0.
\end{equation}
This concludes the proof of Property \itx{B7}{i}.

\vspace{.8cm}

Note that the Weierstrass representations $\Phi_i$ have simple poles and zeroes in $W$. Our next job is to describe a simply connected domain $\Omega$ in $W$ where the above Weierstrass representations determine minimal immersions.

For each $i=1,\ldots,n$, let $D^i$ be an open disk centered at $p_i$ containing $D(p_i,\delta)$. We assume that $D^1,\ldots,D^n$ are pairwise disjoint. Let $\alpha_i\subset D^i\setminus \overline{D(p_i,\delta)}$ be a simple curve connecting $\partial D^i\cap \intc \widehat P$ with $q_i$ and finally let $N_i$ be a small open neighborhood of $\alpha_i\cup\overline{q_ia_i}$ in  $\overline{G_i\cup (D^i\setminus D(p_i,\delta))}$. The domain $\Omega$ is defined as
$$\Omega=\left(\intc \widehat P\setminus\displaystyle\bigcup_{k=1}^n D^k\right)\cup\left(\displaystyle\bigcup_{k=1}^n N_k\right).$$
\begin{figure}[hbtp]
\begin{center}
\includegraphics[width=9cm]{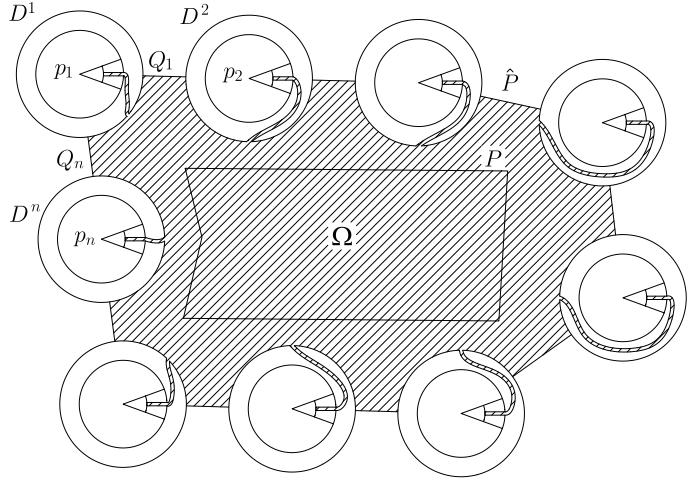}
\caption{The domain $\Omega$.}\label{ques}
\end{center}
\end{figure}

If $D^i$, $\alpha_i$ and $N_i$ are suitably chosen, then  the following properties of $\Omega$ can be guaranteed:
\begin{enumerate}[ ({C}1)]
\item $\overline \Omega$ is simply connected;
\item $\overline{q_ia_i}\subset\overline\Omega$ and $\overline{\intc P}\subset\Omega$;
\item $\overline\Omega$ contains neither the poles, $p_i$, nor the zeroes, $w_i$, of the functions $h_i$, $i=1,\ldots,n$;  
\item $\sup_{z\in \overline{\Omega}} \{\dist_{(\Omega,\eucli)}(0,z)\}<l$. At this point, it is clear the definition of $l$ made in  (\ref{quefrio});
\item $\overline{\Omega}\cap \overline{D(p_i,\delta)}\subset G_i$.
\end{enumerate}

Because of {(C1)} and {(C3)}, the 1-forms of $\Phi^i$ are holomorphic differentials in $\Omega$ without real periods. Therefore, they define conformal minimal immersions in $\Omega$. Even more, it is clear that we can find  a simply connected domain $\Omega'$, with $\overline{\Omega}\subset \Omega'$ and  verifying (C4), where the immersions 
$
X_i:\Omega'\rightarrow\r^3,$
$\displaystyle X_i(z)=\re\int_0^z \Phi^i
$
are well defined, for $i=1, \ldots n.$

Next proposition summarize all the information we will need about the immersions $X_i$.
\begin{proposicion}\label{pancrilver}
For $i=1,\ldots,n$, one has:
\begin{enumerate}[\rm ({D}1$_{i}$)]
\item $\|X_i(z)-X_{i-1}(z)\|<\frac{\epsilon_0}n,\quad \forall z\in\Omega'\setminus D(p_i,\delta)$;
\item $(X_i)_{(3,S_i)}=(X_{i-1})_{(3,S_i)}$;
\item $\|X_n(a_i)-X_n(a_{i+1})\|< 18 \epsilon_0$;
\item $X_n(a_i)\in\r^3\setminus E_{2 \mu}$.
\end{enumerate}
\end{proposicion}

\begin{proof}
The first property is a consequence of  \itx{B6}{i} and {(C4)}. In the same way, \itx{B5}{i} directly implies \itx{D2}{i}. In order to prove \itx{D3}{i} we apply \itx{D1}{k}, $k=1, \ldots, n$, \itx{B7}{i+1}, and (\ref{tarariX}) to obtain 
\begin{multline*}
\|X_n(a_i)-X_n(a_{i+1})\| < \|X_i(a_i)-X_i(q_i)-(X_{i+1}(a_{i+1})-X_{i+1}(q_{i+1}))\|
 +\|X_i(q_i)-X_{i+1}(q_{i+1})\|+2\epsilon_0 <\\
13\epsilon_0+ \| X_0(q_i)-X_0(q_{i+1})\|+4\epsilon_0 \leq 18\epsilon_0.
\end{multline*}

Finally, we prove \itx{D4}{i}.  Using \itx{D1}{i}, \itx{B5}{i},  \itx{B8}{i} (for $z=a_i$), and (\ref{tarariX}), one gets:
\begin{multline*}
\| X_n(a_i)-X(p_i)-3 \mu \nor(X(p_i))\| \leq \|(X_i(a_i)-X_i(q_i))_{(*,S_i)}-3 \mu \nor(X(p_i))\|+|(X_i(a_i)-X_i(q_i))_{(3,S_i)}|+\\ \|X_i(q_i)-X(p_i)\|\leq 
9\epsilon_0.
\end{multline*}
As $X(p_i)+3 \mu \nor(X(p_i)) \in \r^3 \setminus E_{3 \mu}$, then \itx{D4}{i} holds for a small enough  $\epsilon_0$. 
\end{proof}

\subsubsection{The second deformation stage}\label{se:pancrisima}

In this part of the proof, we employ new orthogonal frames. In this case, we take  $T_i=\{w_1^i,w_2^i,w_3^i\}$, $i=1,\ldots,n$, orthonormal bases so that:
\begin{equation}\label{DonJavier}
w_3^i=\nor (X_n(a_i)).
\end{equation}

\begin{figure}[htbp]
        \begin{center}
                \includegraphics[width=0.75\textwidth]{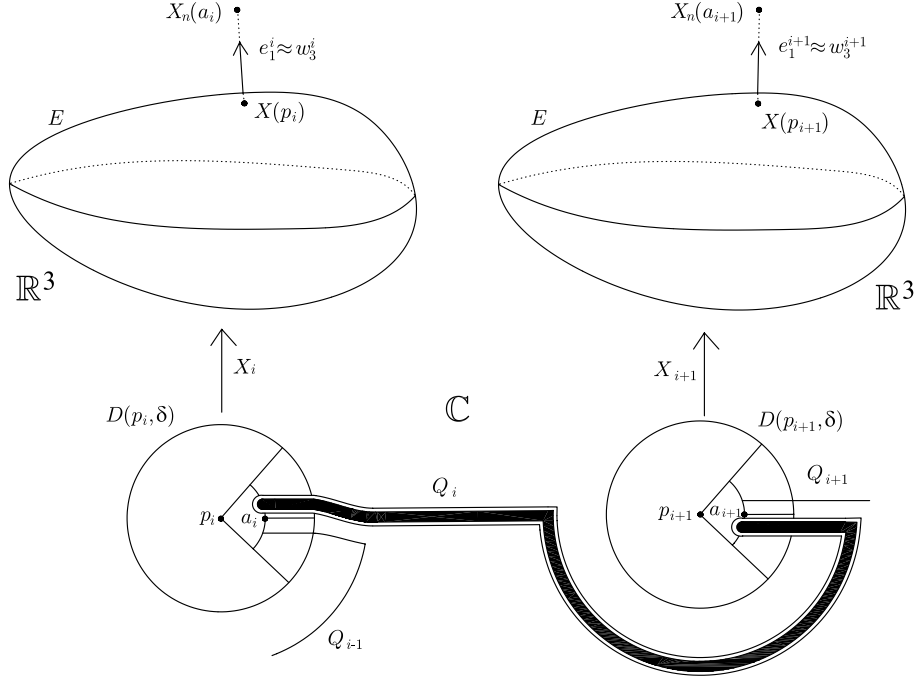}
        \end{center}
        \caption{Immersions $X_i$ and $X_{i+1}$ on the disks $D(p_i,\delta)$ and $D(p_{i+1},\delta)$, respectively.}
        \label{idea2}
\end{figure}

Given $i\in \{1,\ldots,n\}$, we define the curve $Q_i$ as the connected component of $\overline{\partial \Omega\setminus(C_i\cup C_{i+1})}$ that does not cut $C_k$, $k \not\in \{i,i+1\}$ (see Figures \ref{ques} and \ref{idea2}). Note that $Q_i$, $i=1, \ldots,n $ are pairwise disjoint and they satisfy:
\begin{eqnarray}
Q_i&\subset& B^i,\label{sangria}\\
Q_i\cap \overline{D(p_k,\delta)}&=&\emptyset\quad \text{ for }k\not\in\{i,i+1\}.\label{alberto}
\end{eqnarray}
Up to a small perturbation, we can assume that the curves $Q_i$ verify:
\begin{equation}\label{moravia}
f_{(X_n,T_i)}(z)\not=0,\quad \forall z\in Q_i.
\end{equation}

Next step consists of describing neighborhoods of the curves $C_1,\ldots,C_n$, $Q_1,\ldots,Q_n$ that we will use to apply Runge's theorem.

For $i=1\ldots,n$, let $\widehat C_i$ be an open set containing $C_i$ and sufficiently small to fulfill:
\begin{equation}\label{mamanela}
\|X_n(z)-X_n(a_i)\|<3\epsilon_0,\quad \forall z\in \widehat C_i\cap\overline\Omega.
\end{equation}
Notice that the above choice is possible due to Properties \itx{D1}{k} and \itx{B4}{i}. 
Given  $i\in \{ 1,\ldots,n \}$, we define $Q_i^{\xi}=\{z\in\c \: :\: \dist_{(\r^2,\eucli)}(z,Q_i)\leq\xi\},$
where we are assuming that $0<\xi$ is small enough so that:
\begin{enumerate}[ ({E}1)]
\item $Q_i^\xi \subset \Omega'$;
\item $Q_i^\xi\cap Q_j^\xi=\emptyset$, for $i\not=j$;
\item $Q_i^\xi\cap\overline{D(p_k,\delta)}=\emptyset$, for $k\not\in\{i,i+1\}$;
\item $Q_i^{\xi}\subset B^i$;
\item $Q_i^{\xi/2}$ and $\overline{\Omega\setminus Q_i^\xi}$ are simply connected;
\item Given $z \in Q_i$, one has
$|f_{(X_n,T_i)}(z)-f_{(X_n,T_i)}(x)|<\epsilon_1,$ $ \forall x\in B(z,\xi/2),$  where $\epsilon_1=\tfrac14\min_{Q_i}\{|f_{(X_n,T_i)}|\}$;
\item $\sup \{\dist_{\left(\overline{\Omega\setminus Q_i^\xi},\eucli\right)}(0,z)\: : \: z\in\overline{\Omega\setminus Q_i^\xi}\}<l$.
\end{enumerate}
Observe that Properties {(E3)}, {(E4)} and {(E7)} are consequence of (\ref{alberto}), (\ref{sangria}), and {(C4)}, respectively. The other ones are straightforward.

Label $Y_0\stackrel{\text{\tiny def}}=X_n$. We are going to construct a new family of minimal immersions $Y_1,Y_2, \ldots,Y_n$
 with  $Y_i:\Omega'\rightarrow\r^3$ and $Y_i(0)=0$. Associated to these immersions, we will have two families of real parameters $\tau_1,\ldots,\tau_{n-1}$ and $\nu_1,\ldots,\nu_{n-1}$. These families will satisfy the following list of properties, for $i=1,\ldots,n$:
\begin{enumerate}[ ({F}1$_{i}$)]
\item $(Y_i)_{(3,T_i)}=(Y_{i-1})_{(3,T_i)}$;
\item $\|Y_i(z)-Y_{i-1}(z)\|<\frac{\epsilon_0}n$, $\forall z\in\overline{\Omega\setminus Q_i^{\xi}}$;
\item $|f_{(Y_i,T_k)}(z)-f_{(Y_{i-1},T_k)}(z)|<\frac{\epsilon_1}n$, $\forall z\in \overline{\Omega\setminus Q_i^{\xi}}$, for $k=i+1
,\ldots,n$;
\item $ \displaystyle \left(\tfrac{1}{\tau_i}+\tfrac{\nu_i}{\tau_i(\tau_i-\nu_i)}\right)\max_{Q_i^{\xi}}\{|f_{(Y_{i-1},T_i)}g^2_{(Y_{i-1},T_i)}|\}+\nu_i\max_{Q_i^{\xi}}\{|f_{(Y_{i-1},T_i)}|\} <\tfrac2{\xi}$;
\item 
$ \displaystyle \tfrac12\left( \tfrac{\tau_i\xi}{4}\min_{Q_i}\{|f_{(Y_0,T_i)}|\}-1\right)>\text{diam} (E')+1.$
\end{enumerate}
Once again, the sequence $Y_1,\ldots,Y_n$ is defined by recursion. Suppose we have constructed $Y_0,Y_1,\ldots,Y_{i-1}$, for an $i\in\{1,\ldots,n\}$. Hence, we define $Y_i(z)=\re\int_0^z\Phi$ the minimal immersion whose Weierstrass data are given by:
$$\eta_{(Y_i,T_i)}=\eta_{(Y_{i-1},T_i)}l_i,\quad g_{(Y_i,T_i)}=\frac{g_{(Y_{i-1},T_i)}}{l_i},$$
where $l_i:\c\rightarrow\c$ are holomorphic functions satisfying:
\begin{enumerate}[ ({G}1)]
\item $l_i(z)\not=0,\quad\forall z\in\c$;
\item $|l_i(z)-\tau_i|<\nu_i,\quad\forall z\in Q_i^{\xi/2}$;
\item $|l_i(z)-1|<\nu_i,\quad\forall z\in\overline{\Omega\setminus Q_i^{\xi}}$.
\end{enumerate}
Runge's theorem gives us the existence of such functions.
Note that $Y_i:\Omega'\rightarrow \r^3$ has been obtained from $Y_{i-1}:\Omega'\rightarrow\r^3$ by a López-Ros transformation. Notice that $ \phi_{(Y_i,T_k)} \stackrel{\nu_i \to 0}{\longrightarrow} \phi_{(Y_{i-1},T_k)}$ uniformly on $\overline{\Omega \setminus Q^{\xi}_i}$.

Property \itx{F1}{i} trivially holds.
The rest of the properties hold if the constant $\nu_i$ is sufficiently small and $\tau_i$ is large enough. 


\subsubsection{Defining the immersion $Y$}\label{se:ultima}
Define the minimal immersion $Y:\Omega\rightarrow\r^3$ as $Y=Y_n$. We are going to check that $Y$ satisfies the statements of Lemma \ref{lem:propia}.

\paragraph{Statement {(a.2)}:}
Properties {(E4)} and {(A2)} imply that:
$$\overline{\intc P}\subset \Omega\setminus ((\cup_{k=1}^n D(p_k,\delta))\cup(\cup_{k=1}^n Q_k^\xi)).$$
So, we can successively apply \itx{D1}{k} and \itx{F2}{k} to obtain:
\begin{equation}\label{desegundamano}
\|Y(z)-X(z)\|\leq\|Y_n(z)-Y_0(z)\|+\|X_n(z)-X_0(z)\|< 2\epsilon_0<b_1, \quad \forall z\in\overline{\intc P}
\end{equation}
The last inequality occurs if $\epsilon_0$ is small enough.

\paragraph{Statements {(a.1)} and {(a.3)}:}
Consider the following claim:
\begin{claim}\label{risitas}
Every connected curve $\beta$ in $\Omega $ joining the origin $0\in\c$ with $\partial \Omega$ has a point $z'\in\beta$ such that $Y(z')\in \r^3\setminus E'$.
\end{claim}
From (\ref{desegundamano}), it is clear that $Y\left(\overline{\intc P} \right)\subset E'$ (if $\epsilon_0$ is sufficiently small). Then, the existence of a polygon $Q$ verifying {(a.1)} and {(a.3)} is a direct consequence of Claim \ref{risitas}.
\bigskip

\noindent {\em Proof of Claim \ref{risitas}.} 
Let  $\beta\subset\overline\Omega$ be a connected curve with $\beta(0)=0$ and $\beta(1)=z_0\in\partial\Omega$. We will distinguish four possible cases for $z_0$ (see Figures \ref{ques} and \ref{idea2}.)
\begin{itemize}
\item Assume $z_0\in \widehat C_i\cap Q_i^\xi$.

Using \itx{F2}{k} for $k=i+1,\ldots,n$, \itx{F1}{i} and  (\ref{mamanela}), we infer:
$$|(Y_n(z_0)-X_n(a_i))_{(3,T_i)}| \leq 5 \epsilon_0.$$
If we write $T$ as the tangent plane to $\partial E$ at the point $\pro(X_n(a_i))$, then we know that $\dist(p,\partial E)>\dist(p,T)$ for any  $p$ in the halfspace determined by $T$ that does not contain $\partial E$. Taking \itx{D4}{i} into account, and assuming that $\epsilon_0$ is small enough, we have that $Y_n(z_0)$ belongs to the above halfspace, and so: 
\begin{multline}\label{optica}
\dist(Y_n(z_0), \partial E) > \dist (Y_n(z_0),T)=(Y_n(z_0)-\pro(X_n(a_i)))_{(3,T_i)}>\\ (X_n(a_i)-\pro(X_n(a_i)))_{(3,T_i)}-5 \epsilon_0 >  2 \mu-5 \epsilon_0>\mu;
\end{multline} 
here we have used \itx{D4}{i} and (\ref{DonJavier}).
\item Assume $z_0\in \widehat C_i\cap Q_{i-1}^\xi$.

Reasoning as in the above case and using Property \itx{D3}{i-1}, we obtain:
$$|(Y_n(z_0)-X_n(a_{i-1}))_{(3,T_{i-1})}|\leq 23 \epsilon_0.$$
Now, following the arguments of (\ref{optica}), we conclude $Y(z_0)\in\r^3\setminus E'$.
\item The case  $z_0\in \widehat C_i\setminus \cup_{k=1}^n Q_k^{\xi}$ is easier than the previous ones, because (\ref{mamanela}) directly implies 
$\|Y_n(z)-X_n(a_i)\|\leq 4\epsilon_0$, and then we can finish as in the first case.
\item Finally, suppose that $z_0 \in Q_i \setminus \cup_{k=1}^n C_k$.

 Consider $z_1\in\beta\cap\partial D(z_0,\xi/2)$.  For the sake of simplicity, we will write $f^{i-1}$ and  $g^{i-1}$ instead of $f_{(Y_{i-1},T_i)}$ and $g_{(Y_{i-1},T_i)}$, respectively. Furthermore, we will use complex notation to write $a+ {\rm i} b$ instead of $a w_1^i+b w_2^i$. Hence, taking \itx{F2}{k}, $k=i+1, \ldots ,n$, and the definition of $Y_i$ into account one has:
\begin{equation}\nonumber
\begin{split}
&\|Y_n(z_0)-Y_n(z_1)\|\geq \|Y_i(z_0)-Y_i(z_1)\|-2\epsilon_0\geq
 \|(Y_i(z_0)-Y_i(z_1))_{(*,T_i)}\|-2\epsilon_0= \\ &\tfrac12\left|\int_{\overline{z_1z_0}} \overline{f^{i-1}l_idz}-\int_{\overline{z_1z_0}}\frac{f^{i-1}(g^{i-1})^2}{l_i}dz\right|-2\epsilon_0 \geq
\tfrac12\left|\tau_i\int_{\overline{z_1z_0}} \overline{f^{i-1}dz}\right|-\tfrac12\left|\tfrac1{\tau_i}\int_{\overline{z_1z_0}} f^{i-1}(g^{i-1})^2dz\right|-\\
&\tfrac12\left|\int_{\overline{z_1z_0}} \overline{f^{i-1}(l_i-\tau_i)dz}\right|-\tfrac12\left|\int_{\overline{z_1z_0}} f^{i-1}(g^{i-1})^2\left(\frac1{l_i}-\frac1{\tau_i}\right)dz\right|-2\epsilon_0\geq
\end{split}
\end{equation}
using {(G2)}, \itx{F4}{i} and the fact $\longui(\overline{z_1z_0},\eucli)=\xi/2$,
\begin{equation}\nonumber
\begin{split}
&\geq\tfrac{\tau_i}2\left|\int_{\overline{z_1z_0}} \overline{f^{i-1}dz}\right|-\tfrac\xi4\left(\tfrac{1}{\tau_i}\max_{Q_i^{\xi}}\{| f^{i-1}(g^{i-1})^2|\}+\right.\\
&\left.\nu_i\max_{Q_i^{\xi}}\{|f^{i-1}|\}+\frac{\nu_i}{\tau_i(\tau_i-\nu_i)}\max_{Q_i^{\xi}}\{|f^{i-1}(g^{i-1})^2|\}\right)-2\epsilon_0\geq \tfrac12\left(\tau_i\left|\int_{\overline{z_1z_0}} \overline{f^{i-1}dz}\right|-1\right)-2\epsilon_0.
\end{split}
\end{equation}
On the other hand, we make use of \itx{F3}{k}, $k=1,\ldots,i-1$, and {(E6)} to deduce:
\begin{equation}\nonumber
\begin{split}
&\left|\int_{\overline{z_1z_0}} \overline{f^{i-1}(z)dz}\right|\geq\left|f_{(Y_0,T_i)}(z_0)\int_{\overline{z_1z_0}} dz\right|-\\
&\qquad\left|\int_{\overline{z_1z_0}} (f_{(Y_0,T_i)}(z_0)-f_{(Y_0,T_i)}(z))dz\right|-\left|\int_{\overline{z_1z_0}} (f_{(Y_0,T_i)}(z)-f^{i-1}(z))dz\right|\geq \\
&\qquad  \tfrac{\xi}2(|f_{(Y_0,T_i)}(z_0)|-\epsilon_1-\epsilon_1)\geq \tfrac{\xi}2(\min_{Q_i}\{|f_{(Y_0,T_i)}|\}-2\epsilon_1) = \tfrac{\xi}4\min_{Q_i}\{|f_{(Y_0,T_i)}|\}.
\end{split}
\end{equation}
Therefore, by using \itx{F5}{i} we have:
$$\|Y_n(z_0)-Y_n(z_1)\|\geq \tfrac12 \left( \tau_i\tfrac{\xi}{4}\min_{Q_i}\{|f_{(Y_0,T_i)}|\}-1\right)-2\epsilon_0> \text{diam} (E')+1-2\epsilon_0>\text{diam} (E'),
$$
for $\epsilon_0$ small enough.
From the above inequality we conclude that $\beta$ also satisfies the claim in this last case.
\end{itemize}
\hfill{$\square$}

\paragraph{Statement {(a.4)}:} Consider $z\in \intc Q\setminus\intc P$. We will distinguish five possible situations for the complex $z$ (recall that $Q_i^\xi \cap D(p_k,\delta)=\emptyset$, $k \neq i,$ $i+1$)
\begin{itemize}
\item Assume $z\not\in (\cup_{k=1}^n D(p_k,\delta))\cup(\cup_{k=1}^n Q_k^{\xi})$.

With these assumptions Properties \itx{D1}{k} and \itx{F2}{k}), $k=1, \ldots, n$ enable us to conclude that:
$$\|Y_n(z)-X(z)\|\leq 2\epsilon_0<2 b_2,$$
if $\epsilon_0$ is small enough. This fact joint with hypothesis (\ref{eslabon}) give us that $Y(z) \not\in E_{-2 b_2}$.
\item Assume $z\in D(p_i,\delta)\setminus \cup_{k=1}^n Q_k^{\xi}$, for an $i=1,\ldots,n$.
Then, using Properties \itx{D1}{k} $\forall k \neq i$, Properties \itx{F2}{k}, $\forall k$, and (\ref{tarariX}), one has
$$\left< Y_n(z)-X(p_i),e_1^i\right> > \left< X_i(z)-X_i(q_i),e_1^i\right>-4\epsilon_0 >$$
now, using  \itx{B8}{i},
$$>\tfrac12\left(\int_{\overline{q_iz}}\frac{k_idw}{w-p_i}\right)|f_{(\Phi^0,S_i)}(p_i)|-8\epsilon_0>-8 \epsilon_0> -b_2.$$
As a consequence of hypothesis (\ref{eslabon}), we infer $Y_n(z)\not\in E_{-b_2}$. In particular  $Y_n(z)\not\in E_{-2 b_2}$.
\item Assume now $z\in D(p_i,\delta)\cap Q_i^{\xi}$, for an $i=1,\ldots,n$.
As a previous step we need to get an upper bound for $\| e_1^i-w_3^i \|$. In order to do it, notice that:
$$X_n(a_i) \stackrel{\mbox{\scriptsize(D1$_k$)}}{\approx} X_i(a_i) -X_i(q_i)+X_i(q_i)\stackrel{\mbox{\scriptsize (D1$_k$) and (\ref{tarariX})}}{\approx} X_i(a_i)-X_i(q_i)+X(p_i) \stackrel{\mbox{\scriptsize(B8$_i$), (\ref{IA}})}{\approx}$$
$$ 3 \mu e_1^i+(X_i(a_i)-X_i(q_i))_{(3,S_i)} e_3^i+X(p_i) \stackrel{\mbox{\scriptsize (D2$_i$), (D1$_k$) and (\ref{tarariX})}}{\approx} 3 \mu  e_1^i+X(p_i).$$
 Summarizing, it may be checked that:  
$$\|X_n(a_i)- \left(3\mu e_1^i+X(p_i) \right) \| \leq 11\epsilon_0.$$
Therefore, we have:
\begin{equation} \label{milan} \| e_1^i-w_3^i \|=\| \nor(X_n(a_i))-\nor\left(3\mu e_1^i+X(p_i) \right)\| \leq M \|X_n(a_i)- \left(3\mu e_1^i+X(p_i) \right) \| \leq 11M \epsilon_0 ,\end{equation}
where $M$ represents the maximum of $\|d(\nor)\|$ in $\r^3\setminus E$. Note that $M$ does not depend on $\epsilon_0$.

We need another previous inequality:
\begin{multline} \label{fujimori} |\left< Y_i(z)-Y_{i-1}(z),e_1^i\right>|<|\left< Y_i(z)-Y_{i-1}(z),w_3^i\right>+\left< Y_i(z)-Y_{i-1}(z),e_1^i-w_3^i\right>|\leq \\
11M \epsilon_0 \|Y_i(z)\|\|Y_{i-1}(z)\|.\end{multline}
Making use of (\ref{milan}) and (\ref{fujimori}), one obtains:
$$\left< Y_n(z)-X(p_i),e_1^i\right>\geq \left< Y_i(z)-X(p_i),e_1^i\right>-\epsilon_0=$$
$$\left< Y_{i-1}(z)-X(p_i),e_1^i\right>-11M \epsilon_0 \|Y_i(z)\|\|Y_{i-1}(z)\|-\epsilon_0\geq$$
$$\left< X_i(z)-X(p_i),e_1^i\right>-11M \epsilon_0 \|Y_i(z)\|\|Y_{i-1}(z)\|-3\epsilon_0.$$
Reasoning as in the previous case, we conclude:
\begin{equation} \label{mariconazo} \left< Y_n(z)-X(p_i),e_1^i\right> > -b_2-11M \epsilon_0 \|Y_i(z)\|\|Y_{i-1}(z)\|-3\epsilon_0.
\end{equation}
Observe that Statement {\bf (a.3)} and successive applications of \itx{F2}{k}, for $k=i+1, \ldots , n$, give that $Y_i(z) \in E'_{\epsilon_0}$.  Furthermore, notice:
\begin{multline}\|Y_{i-1}(z)-X(q_i)\|\leq \|Y_{i-1}(z)-X_i(z)\|+\| (X_i(z)-X(q_i))_{(*,S_i)}\|+ \\|\left(X_i(z)-X(q_i)\right)_{(3,S_i)}|<3 (\mu +3\epsilon_0),\end{multline}
where have used \itx{F2}{k} and \itx{D1}{k} to get a bound of the first addend, \itx{B8}{i} and \itx{D1}{k} to get a bound of the second addend, and \itx{D2}{i}, \itx{D1}{k} and (\ref{tarariX}) to get a bound of the third one.

Then $\|Y_i(z)\|$ and $\|Y_{i-1}(z)\|$ are bounded in terms of $\epsilon_0$. So, we infer from (\ref{mariconazo}) that  $Y(z)\not\in E_{-2b_2}$, if $\epsilon_0$ is small enough.
\item Suppose  $z\in D(p_{i+1},\delta)\cap Q_i^{\xi}$.

This case is quite similiar to the previous one. Taking into account (\ref{tarariS}) and (\ref{milan}), we deduce $\| e_1^{i+1}-w_3^i \| \leq \| e_1^{i+1}-e_1^i\|+\|e_1^i-w_3^i\|<\epsilon_0(11 M+1)$.
As in (\ref{fujimori}), we obtain:
$$|\left< Y_i(z)-Y_{i-1}(z),e_1^{i+1}\right>|<|\left< Y_i(z)-Y_{i-1}(z),w_3^i\right>+\left< Y_i(z)-Y_{i-1}(z),e_1^{i+1}-w_3^i\right>|\leq$$
$$ \epsilon_0 (11 M+1)\|Y_i(z)\|\|Y_{i-1}(z)\|$$
Using these inequalities as in the former case, we deduce $Y(z)\not\in E_{-2b_2}$.
 
\item Finally, assume $z\in Q_i^{\xi}\setminus\cup_{k=1}^n D(p_k,\delta)$.

If we apply the inequalities that we got in the previous cases, then we have:
$$\left< Y_n(z)-X(p_i),e_1^i\right>\geq \left< Y_i(z)-X(p_i),e_1^i\right>-\epsilon_0=$$
$$\left< Y_{i-1}(z)-X(p_i),e_1^i\right>-11 M \epsilon_0 \|Y_i(z)\|\|Y_{i-1}(z)\|-\epsilon_0\geq$$
$$\left< X(z)-X(p_i),e_1^i\right>-11 M \epsilon_0\|Y_i(z)\|\|Y_{i-1}(z)\|-3\epsilon_0\geq$$
using (\ref{tarariX}):
$$\geq \epsilon_0-11 M \epsilon_0 \|Y_i(z)\|\|Y_{i-1}(z)\|-3\epsilon_0\geq-2b_2.$$
\end{itemize}

This concludes the proof of Statement (a.4) and completes the proof of Lemma \ref{lem:propia}.


\subsection{Proof of Lemma \ref{lem:nadi}}

The proof of this lemma is a  refinement of Nadirashvili's techniques (see \cite{nadi}.) However, we must adapt his ideas about balls to the setting of general convex bodies of space. Hence, our deformation increases the intrinsic diameter, but it controls the extrinsic geometry in $\r^3$, in such a way that the perturbed surface is included in a prescribed convex set. Lemma \ref{lem:nadi}  we will be the key in the proof of Theorem \ref{th:mari}. 

As in the former section, we denote $\nor$ and $\pro$ as the extended Gauss map associated to $E$ and the normal projection to $\partial E$, respectively. We also need some notation about angles between vectors of $\r^3$. Given $a,b\in\r^3$ and $\theta\in [ 0, \pi[$, we label the angle formed by $a$ and $b$ as $\angle(a,b) \in [ 0, \pi[$  and write $\mbox{Cone}(a,\theta)= \left\{p \in \r^3 \; | \; \angle (p,a) < \theta \right\}.$

Consider $P$, the polygon given in the statement of Lemma \ref{lem:nadi}. 
Our first step is to describe a labyrinth on $\intc P$, which depends on $P$ and a positive integer $N$. Let $\ell$  be the number of sides of $P$. From now on $N$ will be a  positive multiple of $\ell$.

\begin{remark}
Throughout the proof of the lemma a set of real positive constants  depending on $X$, $P$,  $\varepsilon$, $\sigma$, $a$, and $b$ will appear. The symbol `$\cte$' will denote these different constants. It is important to note that the choice of these constants does not depend on the integer $N$.
\end{remark}

\begin{figure}[htb] 
\begin{center}
\includegraphics[width=7.8cm]{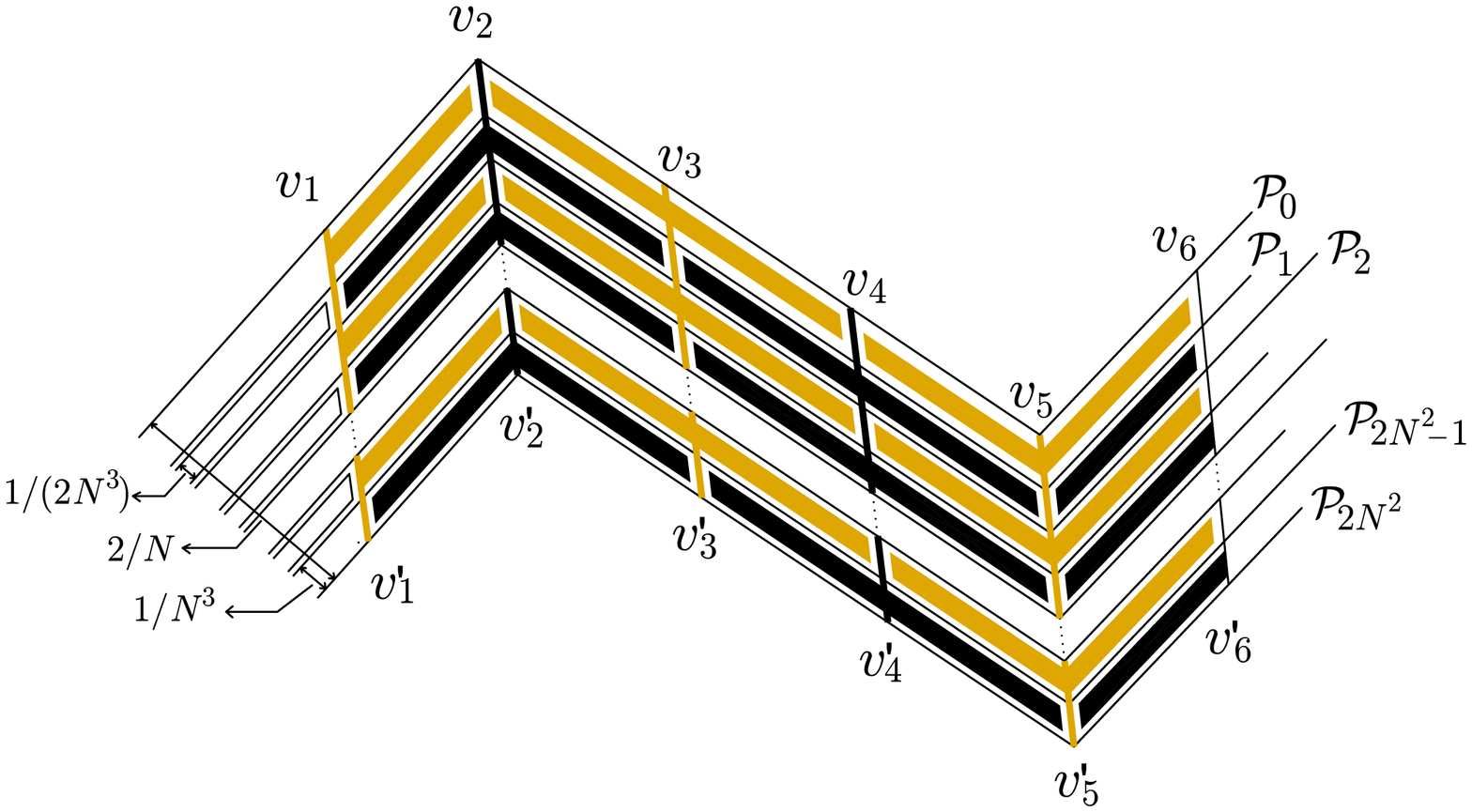} \hspace{.8cm}
\includegraphics[width=6cm]{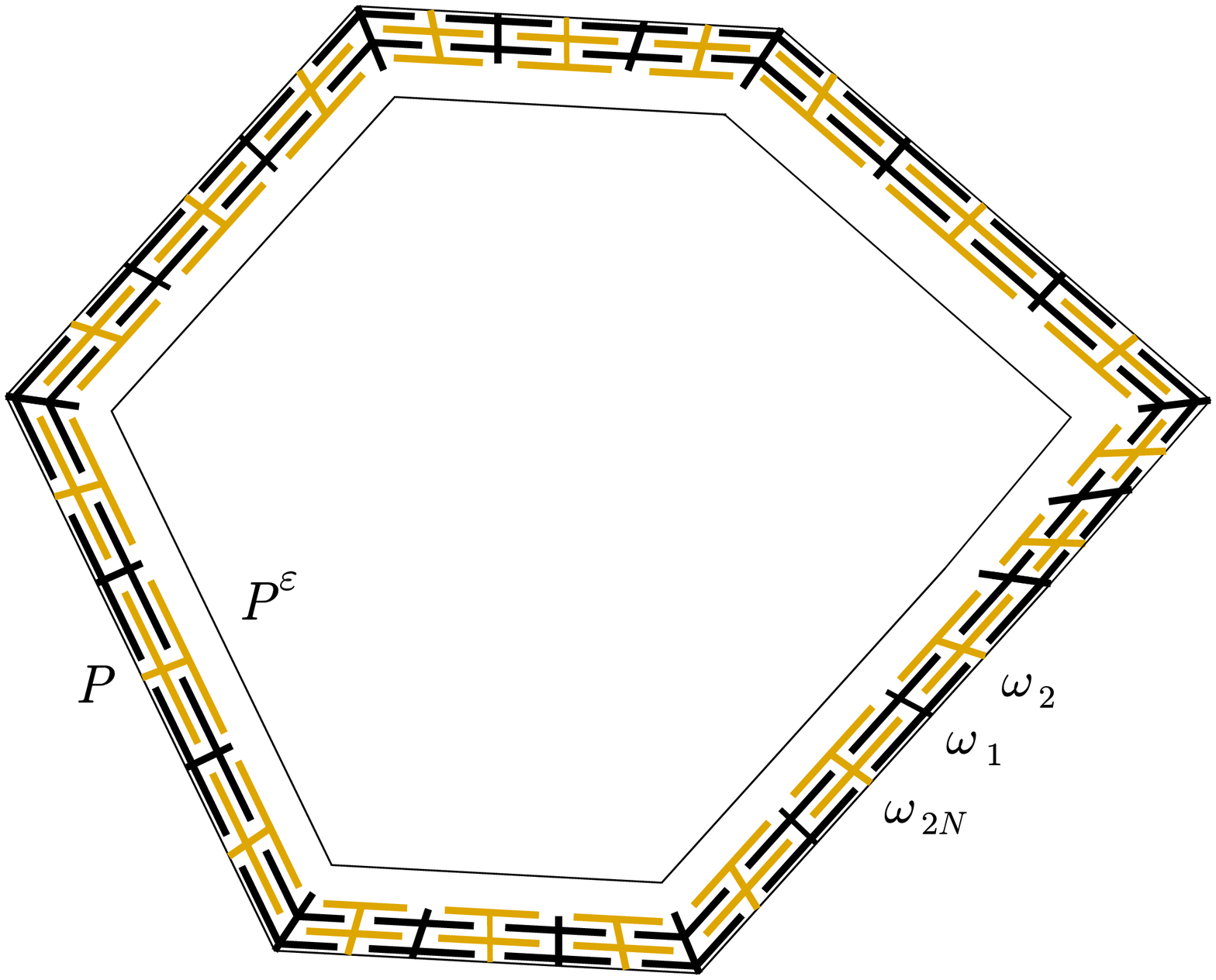}
\caption{Distribution of the sets $\omega_i$. \label{fig:dos}}
\end{center}
\end{figure}

Let $\zeta_0>0$ small enough so that $P^{\zeta_0}$ is well defined and $\overline{\intc(P^\varepsilon)} \subset \intc(P^{\zeta_0})$. From now on, we will only consider $N \in \n$ such that $2/N < \zeta_0$.
Let $v_1, \ldots ,v_{2N}$ be a set of points in the polygon $P$ (containing the vertices of $P$) that divide each side of $P$ into $\frac{2N}{\ell}$ equal parts. We can transfer this partition to the polygon $P^{2/N}$: $v_1^\prime, \ldots, v_{2N}^\prime$. We define the following sets:
\begin{itemize}
\item $L_i=$ the segment that joins $v_i$ and $v_i^\prime$, $i=1, \ldots 2 N$;
\item ${\cal P}_i=P^{i/N^3}, \; i=0, \ldots 2N^2$;
\item $\mathcal{E}= \bigcup_{i=0}^{N^2-1} \overline{\intc({\cal P}_{2i}) \setminus \intc({\cal P}_{2i+1})}$ and $\widetilde{\mathcal{E}}= \bigcup_{i=1}^{N^2} \overline{\intc({\cal P}_{2i-1}) \setminus \intc({\cal P}_{2i})}$;
\item $\mathcal R= \bigcup_{i=0}^{2N^2} {\cal P}_i$;
\item $\mathcal{B}= \bigcup_{i=1}^N L_{2i}$ and $\widetilde{\mathcal{B}}= \bigcup_{i=0}^{N-1} L_{2i+1}$;
\item $L=\mathcal{B} \cap \mathcal{E}$, $\widetilde{L}=\widetilde{\mathcal{B}} \cap \widetilde{\mathcal{E}}$, and $H=\mathcal R \cup L \cup \widetilde{L}$;
\item $\Xi_N= \{ z \in \intc({\cal P}_0) \setminus \intc({\cal P}_{2N^2}) : \dist_{(\c,\langle \cdot,\cdot \rangle)}(z,H) \geq \frac{1}{4N^3}\}$.
\end{itemize}
We define $\omega_i$ as the union of the segment $L_i$ and those connected components of $\Xi_N$ that have nonempty intersection with $L_i$ for $i=1, \ldots,2 N$. Finally, we label $\varpi_i= \{ z \in \c \: : \: \dist_{(\c,\langle \cdot,\cdot \rangle)}(z,\omega_i)< \delta(N) \}$, where   $\delta(N)>0$ is chosen in such a way that the sets $\overline{\varpi_i}$ ($i=1, \ldots ,2 N$) are pairwise disjoint (see Figure \ref{fig:dos}).

The shape of the labyrinth formed by the sets $\omega_i$ guarantee
the following claims hold if $N$ is large enough:
\begin{description} 
\item [{\rm\em Claim A.}] The Euclidean diameter of $\varpi_i$ is less than $\frac{\cte}{N}$. 
\item [{\rm\em Claim B.}] $\nor(X(\varpi_i))\subset \text{Cone} (q,\frac 1{\sqrt N})$ for a suitable $q\in\s^2$. This is a consequence of Claim A.
\item [{\rm\em Claim C.}] If $\lambda^2 \left< \cdot , \cdot \right>$ is a conformal metric on $\overline{\intc P}$ and satisfies:
$$\lambda \geq\begin{cases}
 c & \text{in } \intc P,\\
 cN^4 &\text{in } \Xi_N,
\end{cases}$$ for $c \in \r^+$,
and if $\alpha$ is a curve in $\overline{\intc P}$ connecting $P^\varepsilon$ and $P$, then
$\longui (\alpha,\lambda \left<\cdot,\cdot\right>)>\cte \, c \, N,$ where $\cte$ does not depend on $c$. \newline
Claim C is a consequence of the fact that a curve $\alpha$, that does not go through the connected components of $\Xi_N$, must have a large Euclidean length.
\end{description}

As we have done several times in the previous section, we deform $X$ in $2N$ steps, by using López-Ros transformation and Runge's theorem. Then, we construct (for a sufficiently large $N$)  a sequence of $2N$ conformal minimal immersions defined on $\overline{\intc P}$: $F_0=X, F_1, \ldots ,F_{2N}$, which satisfies: 
\begin{enumerate}[\rm ({H}1$_{i}$)]
\item $F_i(z)= \re \left( \int_0^z \phi^i(w)dw \right)$;
\item $\|\phi^i(z)-\phi^{i-1}(z) \| \leq 1/N^2$ for all $z \in \overline{\intc P} \setminus \varpi_i$;
\item $\|\phi^i(z) \| \geq N^{7/2}$ for all $z \in \omega_i$;
\item $\| \phi^i(z) \| \geq \frac \cte{\sqrt{N}}$ for all $z \in \varpi_i$;
\item $\dist_{\s^2}(G_i(z), G_{i-1}(z))<\frac{1}{N^2}$ for all $z \in \intc P \setminus \varpi_i$, where $\dist_{\s^2}$ is the intrinsic distance in $\s^2$ and $G_i$ represents the Gauss map of the immersion $F_i$;
\item there exists a set $S_i=\{ e_1, e_2, e_3 \}$ of orthonormal frame in $\r^3$ such that:
\begin{enumerate}[\rm ({H6}.1$_{i}$)]
\item For any $z \in \overline{\varpi_i}$,  we have $\| \nor(X(z))_{(*,S_i)} \| <\frac {\cte}{\sqrt{N}}$;
\item $(F_i(z))_{(3,S_i)}=(F_{i-1}(z))_{(3,S_i)}$ for all $z \in \overline{\intc P}$;
\end{enumerate}
\item $\| F_i(z)-F_{i-1}(z) \| \leq \frac{\cte}{N^2}$, $\forall z \in \intc P \setminus \varpi_i$.
\end{enumerate}


The sequence $F_0,F_1,\ldots,F_{2N}$ is constructed in a recursive way. 
Suppose that we have $F_0, \ldots, F_{j-1}$ verifying the claims (H1$_i)$, $ \ldots,$ (H7$_i$), $i=1, \ldots,$ $j-1$.
In order to construct $F_j$, first note that (for a  large enough $N$) the following statements hold:
\begin{enumerate}[({J}1)]
\item $\|\phi^{j-1}\| \leq \cte$ in $\intc P \setminus \bigcup_{k=1}^{j-1} \varpi_k$. \newline We easily deduce this property from (H2$_l$) for $l=1, \ldots, j-1$. 
\item $\|\phi^{j-1}\| \geq \cte$ in $\intc P \setminus \bigcup_{k=1}^{j-1} \varpi_k$. \newline To obtain this property, it suffices to apply (H2$_l$) for $l=1, \ldots, j-1$ once again.
\item The diameter in $\r^3$ of $F_{j-1}(\varpi_j)$ is less than $\frac 1{\sqrt{N}}$. \newline This is a consequence of (J1), Claim A, and (\ref{eq:metric}).
\item The diameter in $\s^2$ of $G_{j-1}(\varpi_j)$ is less than $\frac 1{\sqrt{N}}$. In particular $G_{j-1}(\varpi_j) \subset \hbox{Cone}\left(g , \frac{1}{\sqrt{N}}\right)$, for some $g \in G_{j-1}(\varpi_j)$. \newline From Claim A, the diameter of $G_0(\varpi_j)$ is bounded. Then (J4) holds after successive applications of (H5$_l$).
\item There exists an orthogonal frame $S_j=\{e_1, e_2, e_3 \}$ in $\r^3$, where: 
\begin{enumerate}[({J5.}1)] 
\item $\angle \left(e_3, \nor(X(z))\right) \leq \tfrac{\cte}{\sqrt{N}}$ for all $z \in \varpi_j$;
\item $\angle (\pm e_3,G_{j-1}(z)) \geq \tfrac{\cte}{\sqrt{N}}$ for all $z \in \varpi_j$.
\end{enumerate}
The proof of (J5) is slightly more complicated. Let  $C=\hbox{Cone}\left(g,\frac{2}{\sqrt N}\right)$  where $g$ is given by Property (J4). To obtain (J5.2) it suffices to take $e_3$ in  $\s^2 \setminus H$, where
$H=C \bigcup (-C).$
On the other hand, in order to verify (J5.1), the vector $e_3$ must be chosen as follows:
{\samepage
\begin{itemize}
\item If $(\s^2\setminus H)\cap \nor (X(\varpi_j))\not=\emptyset$, then we take $e_3\in(\s^2\setminus H)\cap \nor (X(\varpi_j))$; 
\item If $(\s^2\setminus H)\cap \nor (X(\varpi_j))=\emptyset$, then we take $e_3\in \s^2 \setminus H$ satisfying $\angle(e_3,q')<\frac{2}{\sqrt N}$ for some $q'\in \nor (X(\varpi_j))$.
\end{itemize}
}
It is straightforward to check that this choice of $e_3$ guarantees (J5).
\end{enumerate}
We shall now construct $F_j$.  Let $(g^{j-1},\phi_3^{j-1})$ be the Weierstrass data of the immersion $F_{j-1}$ in the coordinate system $S_j$. Let $h_\alpha:\c \rightarrow \c \setminus\{0\}$ be a holomorphic function verifying:
\begin{itemize}
\item $|h_\alpha(z)-1|<1/\alpha$, $\forall z \in \overline{\intc P} \setminus \varpi_j$;
\item $|h_\alpha(z)-1/\alpha|<1/\alpha$, $\forall z \in \omega_j,$ 
\end{itemize}
where $\alpha$ is a sufficiently large positive number.
The existence of such a function is given by Runge's theorem.

We define ${\phi}^j_3=\phi_3^{j-1}$ and $g^j=g^{j-1}/h_\alpha$. Let $F_j$ be  the minimal immersion  induced by the aforementioned Weierstrass data in the frame $S_j$, $F_{j}(z)= \re \left( \int_0^z \phi^{j}(w) \, dw \right).$

$F_j$ must now verify the properties (H1$_j$),$ \ldots $,(H7$_j$). (Note that they do not depend on changes of coordinates in $\r^3$). Claim (H1$_j$) directly follows from the definition.

Note that $h_\alpha \to 1$ (resp. $h_\alpha \to \infty$) uniformly on $\overline{\intc P} \setminus \varpi_j$ (resp. on $\omega_j$), as $\alpha \to \infty$. Then 
(H2$_j$), (H3$_j$), (H5$_j$), and (H7$_j$) easily hold for a large enough $\alpha$  (in terms of $N$.)

(H4$_j$) is verified using (J5.2) which gives:
$$\frac{\sin \left(\frac{\cte}{\sqrt{N}}\right)}{1+\cos\left(\frac{\cte}{\sqrt{N}}\right)} \leq |g^{j-1}| \leq \frac{\sin \left(\frac{\cte}{\sqrt{N}}\right)}{1-\cos\left(\frac{\cte}{\sqrt{N}}\right)} \qquad \hbox{in } \varpi_j,$$
and so, taking (J2) into account one has (if $N$ is large enough):
$$\| \phi^j \| \geq | \phi^j_3|=| \phi^{j-1}_3| \geq \sqrt2\| \phi^{j-1} \| \frac{|g^{j-1}|}{1+|g^{j-1}|^2}\geq \cte \cdot \sin\left(\tfrac{\cte}{\sqrt{N}}\right) \geq \tfrac \cte{\sqrt{N}} \qquad \hbox{in } \varpi_j.$$

Using (J5.1), we get (H6.1$_j$). To obtain (H6.2$_j$), we use $\phi^{j-1}_3=\phi_3^j$ in the frame $S_j$.

Hence, we have constructed the immersions $F_0,F_1, \ldots ,F_{2N}$ verifying claims (H1$_j$),$\ldots$,(H7$_j$) for $j=1, \ldots , 2N$. 

In the following claim we establish some properties of the immersion $F_{2N}$ that will be essential for our purpose.
\begin{claim}
 If $N$ is large enough, then $F_{2N}$ verifies that:
\begin{enumerate}[\rm ({K}1)]
\item $2\sigma< \dist_{(\overline{\intc P},\metri{F_{2N}})}(P,P^\varepsilon)$;
\item $\|F_{2N}(z)-X(z) \| \leq \frac{\cte}{N}, \forall z\in  \intc P \setminus \bigcup_{j=1}^{2N} \varpi_j$;
\item there is a polygon $\widetilde{P}$ such that
\begin{enumerate}[\rm ({K3}.1)]
\item $\overline{\intc P^\varepsilon} \subset \intc{\widetilde{P}} \subset \overline{\intc{\widetilde{P}}} \subset \intc{P}$;
\item  $\sigma<\dist_{(\overline{\intc P},\metri{F_{2N}})}(z,P^\varepsilon)<2\sigma$, $\forall z \in \widetilde{P}$;
\item $F_{2N}(\overline{\intc \widetilde P}) \subset E_{\mu'/2}$, where $\mu'=\dist (\partial E, \partial E')$;
\item $F_{2N}(\overline{\intc P\setminus \intc P^\varepsilon}) \subset \r^3 \setminus E_{-2 a}$.
\end{enumerate}
\end{enumerate}
\end{claim}
\begin{proof}
To get Property (K1) we have to apply {\em Claim C}, taking into account properties (H2$_k$), (H3$_k$), and (H4$_k$) (see Proposition 2 in \cite{cili} for details). Property (K2) is deduced from a successive use of (H7$_k$), for $k=1, \ldots, 2N$.

We now  construct the polygon $\widetilde P$ of Property (K3). Let
$$\mathcal{A}=\{ z \in \intc P \setminus \intc P^\varepsilon \: : \:\sigma<\dist_{(\overline{\intc P},\metri{F_{2N}})}(z,P^\varepsilon)<2\sigma \}.$$
$\mathcal{A}$ is a nonempty open subset of $\intc P\setminus\intc P^\varepsilon$.
Observe that the polygons $P$  and $P^\varepsilon$ are contained in
different path components of $\c \setminus \mathcal{A}$. Then, we can draw a polygon $\widetilde{P}$ on $\mathcal{A}$ such that $\overline{\intc P^\varepsilon} \subset \intc \widetilde{P}\subset \overline {\intc \widetilde P}\subset \intc P$.
Thus, (K3.1) and (K3.2) trivially hold.

Item (K3.3) is checked in the following paragraphs.
Thanks to the convex hull property (see \cite{osserman2}), we only need to prove that 
$F_{2N}(\widetilde{P}) \subset E_{\mu'/2}.$

Let $\eta \in \widetilde{P}$. If $\eta \in \intc P \setminus \bigcup_{j=1}^{2N} \varpi_j$, then we have $ \| F_{2N} (\eta)-X(\eta)\|  \leq \frac{\cte}{N} $, and we know that $X(\eta) \in E$. Therefore, $F_{2 N}(\eta) \in E_{\frac{\cte}{N}} \subset E_{\mu'/2}$, if $N$ is large enough.

On the other hand, if $\eta \in \varpi_j$ for $j\in \{1, \ldots , 2N\}$ then to check that $F_{2N}(\eta) \in E_{\mu'/2}$  is slightly more complicated. From (K3.2), it is possible to find a curve $\gamma :[0,1] \rightarrow \intc P$ such that $\gamma(0) \in P^\varepsilon, \gamma(1)=\eta$, and $\longui(\gamma,\metri{F_{2N}}) \leq 2\sigma$. We define:
$$\overline{t}=\sup \{ t \in [0,1] \: : \: \gamma(t) \in \partial \varpi_j \}, \qquad \overline{\eta}=\gamma(\overline{t}).$$
Notice that the previous supremum exists because  $\varpi_j \subset \intc(P) \setminus \overline{\intc(P^\varepsilon)}$ (for a large enough $N$.) 
To continue, we need to demonstrate that
\begin{equation} \label{oeta}
\|F_j(\overline{\eta})-F_j(\eta)\| \leq \frac{\cte}{N}+2\sigma.
\end{equation}
Indeed, taking Properties (H7$_k$) into account, one has:
\begin{multline}\nonumber
\| F_j(\overline{\eta})-F_j(\eta)\| \leq \| F_j(\overline{\eta})-F_{2N}(\overline{\eta}) \| + \|F_{2N}(\overline{\eta})-F_{2N}(\eta)\| +\|F_{2N}(\eta)-F_j(\eta)\| \leq \\
\leq \frac{\cte}N+\longui(\gamma,\metri{F_{2N}})+\frac{\cte}N \leq \frac{\cte}N+2\sigma.
\end{multline}
Now, define $\Sigma(\eta)= \pro(X(\eta))-\frac1{\kappa_2(\partial E)} \nor(X(\eta)),$ and $\s^2(\eta)$ the sphere centered at $\Sigma(\eta)$ of radius $\frac1{\kappa_2(\partial E)}$. Notice that $\partial E$ is contained in the exterior of $\s^2(\eta)$. 
\begin{figure}[h]
        \begin{center}
                \includegraphics[width=0.50\textwidth]{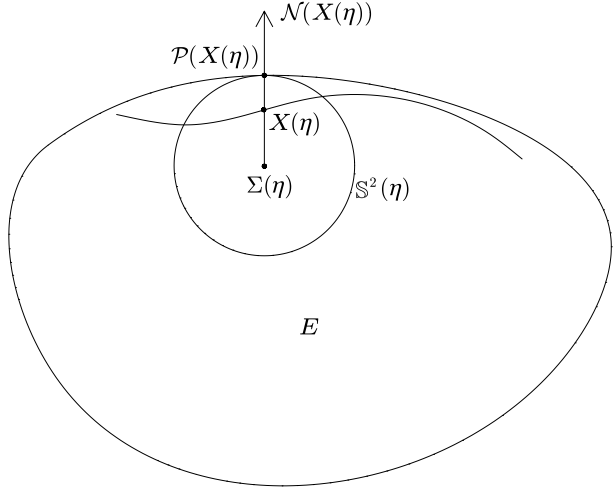}
        \end{center}
        \caption{The sphere $\s^2(\eta)$.}
        \label{fig:mingo}
\end{figure}

In what follows, we will assume that $F_{2N}(\eta)\not\in E$, if not (K3.3) trivially holds. Our next job is to find an upper bound for the distance $\dist(F_{2N}(\eta),\partial E)$. 
$$\dist(F_{2N}(\eta),\partial E)\leq \dist(F_{2N}(\eta),\s^2(\eta))\leq \dist(F_{j}(\eta),\s^2(\eta))+\frac{\cte}N=$$
$$\|F_j(\eta)-\Sigma(\eta)\|-\frac1{\kappa_2(\partial E)}+\frac{\cte}N=\sqrt{(F_j(\eta)-\Sigma(\eta))_{(3,S_j)}^2+\|(F_j(\eta)-\Sigma(\eta))_{(*,S_j)}\|^2}-\frac1{\kappa_2(\partial E)}+\frac{\cte}N.$$
We separately bound each addend in the root. For the third coordinate, we use Hypothesis 1), (H6.2$_j$), and (H7$_k$)
$$|(F_j(\eta)-\Sigma(\eta))_3|=|(F_{j-1}(\eta)-\Sigma(\eta))_3| \leq |(F_{j-1}(\eta)-X(\eta))_3|+|(X(\eta)-\Sigma(\eta))_3| \leq \frac{\cte}{N}+a+\frac1{\kappa_2(\partial E)}.$$

Regarding the first two coordinates, we can apply inequality (\ref{oeta}), Properties (H7$_k$), (J3) and (H6.1$_j$) (taking Hypothesis 1) into account) and so one has 
\begin{multline*} \|(F_j(\eta)-\Sigma(\eta))_* \| \leq \| (F_j(\eta) - F_j(\overline{\eta}))_* \|
+\| (F_j(\overline{\eta}) - F_{j-1}(\overline{\eta}))_*\| +\| (F_{j-1}(\overline{\eta}) -F_{j-1}(\eta))_*\|+ \\ 
\|(F_{j-1}(\eta)-X(\eta))_*\| 
+ \|(X(\eta)-\pro(X(\eta))_*\|  +\|\frac 1{\kappa_2(\partial E)} \nor(X(\eta))_*\|\leq \\
 \frac{\cte}{N}+2\sigma+\frac{\cte}{N^2}+\frac{1}{\sqrt N}+\frac{\cte}{N}+a+ \frac1{\kappa_2(\partial E)} \frac{\cte}{\sqrt{N}}\leq  2 \sigma+a+\frac{\cte}{\sqrt{N}}.\end{multline*}
Using Pythagoras' theorem and Hypothesis 2), one obtains:
$$\dist(F_{2N}(\eta),\partial E) \leq \sqrt{\left(\frac{\cte}{N}+a+\frac1{\kappa_2(\partial E)}\right)^2+\left(2 \sigma+a+\frac{\cte}{\sqrt{N}}\right)^2}-\frac1{\kappa_2(\partial E)}+\frac{\cte}N<$$ $$\sqrt{\left(a+\frac1{\kappa_2(\partial E)}\right)^2+\left(2 \sigma+a\right)^2}-\frac1{\kappa_2(\partial E)}+\varepsilon=\frac{\mu'}2$$
for a large enough $N$.

Finally, we check (K3.4). Let $z\in \overline{\intc P\setminus \intc P^\varepsilon}$. If $z\not \in\bigcup_{k=1}^{2N}\varpi_k$, then we deduce (K3.4) from  Hypothesis 1) in the lemma and item (K2).
If $z\in\varpi_j$, we proceed as follows. We consider the frame $S_j=\{ e_1,e_2,e_3 \}$, and $T$ the tangent plane of $\partial E_{-2 a}$ at $\pro_{E_{-2 a}}(X(z))$. Then, a lower bound of the oriented distance from $F_{2N}(z)$ to $T$ can be estimated in the following way:
$$ \langle F_{2N}(z)-\pro_{E_{-2 a}}(X(z)),\nor(X(z)) \rangle = \langle F_{2N}(z)-\pro_{E_{-2 a}}(X(z)),e_3 \rangle+\langle F_{2N}(z)-\pro_{E_{-2 a}}(X(z)),\nor(X(z))-e_3 \rangle \geq$$
applying (H6.2$_j$) and (H7$_k$), $k=j+1, \ldots, 2N,$ one has
$$ \geq \langle F_{j-1}(z)-\pro_{E_{-2 a}}(X(z)),e_3 \rangle-\| F_{2N}(z)-\pro_{E_{-2 a}}(X(z))\| \cdot \|\nor(X(z))-e_3 \|-\frac{\cte}{N}>$$
now, we use (K3.3), (H6.1$_j$) and (H7$_k$) to obtain
\begin{multline*}\geq \langle F_{j-1}(z)-\pro_{E_{-2 a}}(X(z)),e_3 \rangle-\frac{\cte}{\sqrt N}
\geq\langle X(z)-\pro_{E_{-2 a}}(X(z)),e_3 \rangle-\frac{\cte}{\sqrt N} \geq \\ \langle X(z)-\pro_{E_{-2 a}}(X(z)),\nor(X(z)) \rangle+\langle X(z)-\pro_{E_{-2 a}}(X(z)),e_3-\nor(X(z)) \rangle-\frac{\cte}{\sqrt{N}}.\end{multline*}
Finally, Hypothesis 1), and Property (H6.1$_j$) imply
$$\langle F_{2N}(z)-\pro_{E_{-2 a}}(X(z)),\nor(X(z)) \rangle \geq a-\frac{\cte}{\sqrt{N}}>0.$$
Thus, $\dist(F_{2 N}(z), \partial(E_{-2 a}))>0,$ if $N$ is large enough. In particular, $F_{2 N}(z) \not\in E_{-2 a}$.
This completes the proof of the claim. \end{proof}

The immersion $F_{2N}$ and the polygon $\widetilde P$ verify all the properties of Lemma \ref{lem:nadi} except for (b.3). In order to get this property, we modify $F_{2N}$ by using Lemma \ref{lem:propia}.

Hence, we apply Lemma \ref{lem:propia} to the following data:
$$\hat X=F_{2N},\quad \hat P=\widetilde P,\quad \hat E=E_{-2 a},\quad \hat E'=E',\quad\hat b_2=b,$$
 and $\hat O\subset \intc P$ is an open neighborhood of $\overline
 {\intc \widetilde P}$ which is small enough to verify
 (\ref{eslabon}). We can choose $\hat O$ because of (K3.3) and (K3.4).  Then, we consider a polygon $Q$ and a minimal immersion $Y:\overline{\intc Q}\rightarrow \r^3$ given by applying Lemma \ref{lem:propia} to the above data and a suitable $\hat b_1$. It is easy to check that $Y$ and $Q$ verify all the claims in Lemma \ref{lem:nadi}. 
Claim (b.2) is a direct consequence of (K3.2) and the fact that $Y$ converges to $F_{2N}$, uniformly on $\overline{ \intc \widetilde P}$, as $\hat b_1\rightarrow 0$.
Claim (b.4) is deduced as follows. Item (a.4) in Lemma \ref{lem:propia} implies $Y(\intc Q\setminus \intc \widetilde P)\subset \r^3\setminus E_{-2(a+b)}$. On the other hand, using (K3.4) and item (a.2) in Lemma \ref{lem:propia}, we deduce $Y(\intc \widetilde P\setminus \intc P^\varepsilon)\subset \r^3\setminus E_{-2a-
\hat b_1}$. This completes the proof of Lemma \ref{lem:nadi}.


\section{Proof of the main Theorems}\label{sec:theorem}
In this section we state and prove the principal results of this
article. First we demonstrate the existence of properly immersed
minimal disks in bounded regular convex domains, i.e., domains whose boundary is a compact analytic surface of $\r^3$ (Theorem \ref{th:mari}.) Later, we make use of a classical result by Minkowski (Theorem \ref{th:minko}) to construct complete proper minimal disks in any convex domain (Theorem \ref{th:refinitivo} and Corollary \ref{co:esperanza}).
\begin{teorema} \label{th:mari}
Let $D $ and $D'$ be two bounded convex regular domains satisfying $\overline{D} \subset D'$. Let $\varphi:\overline{\intc(\Pi)} \rightarrow   \r^3$ be a conformal minimal immersion, where $\Pi \subset \d$ is a polygon. Assume that $\varphi(\Pi)\subset D \setminus D_{-d}$, where $d$ is a positive constant such that:
$$\sqrt{\left(d+\frac1{\kappa_2(\partial D)}\right)^2+d^2}-\frac1{\kappa_2(\partial D)}<\frac{\dist(\partial D,\partial D')}2.$$
Then, for any $\delta>0$, there exist a simply connected domain $\Omega \subset \d$, with $\overline{\intc(\Pi)} \subset \Omega$, and a complete proper minimal immersion $\widetilde \varphi: \Omega \rightarrow D'$ such that: 
\begin{enumerate}[\rm (a)]
\item $\|\widetilde \varphi-\varphi\|< \delta , \quad \text{in } \intc(\Pi);$
\item $\widetilde \varphi(\Omega \setminus \intc \Pi) \subset D' \setminus D_{-2d-\delta}.$
\end{enumerate}
\end{teorema}

\begin{proof}
First of all, we define a sequence $\{E^n \}$ of bounded convex regular domains in the following way. Consider $\nu>0$ small enough to verify that $D'_{-\widetilde\nu}$ exists, $\overline{D}\subset D'_{-\widetilde\nu}$, and
$$\sqrt{\left(d+\frac1{\kappa_2(\partial D)}\right)^2+d^2}-\frac1{\kappa_2(\partial D)}<\frac{\dist(\partial D,\partial D'_{-\widetilde \nu})}2, $$
where $\widetilde\nu=\sum_{k=2}^\infty \nu/k^2$. Then, we define $E^1\df D$ and $E^n\df D'_{-\sum_{k=n}^\infty \nu/k^2}$, $n\geq 2$. We also take a decreasing sequence of positive reals $\{b_n\}$ with  $b_1=d$, $b_n<1/n^3$ for $n>1$, and which is related to $\{E^n\}$ in the following way
$$\sqrt{\left(b_{n}+\frac1{\kappa_2(\partial E^{n})}\right)^2+(b_{n})^2}-\frac1{\kappa_2(\partial E^{n})}<\frac{\mu'_{n}}2, \qquad \mbox{where $\mu'_n=\dist(\partial E^n, \partial E^{n+1}).$}$$

Next, we use Lemma \ref{lem:nadi} to construct a sequence 
$$\chi_n=(X_n:\overline{\intc P_n} \rightarrow \r^3,P_n, \varepsilon_n,\xi_n,\sigma_n),$$
where $X_n$ are conformal minimal immersions with $X_n(0)=0$, $P_n$ are polygons,  and $ \{     \varepsilon_n\}$, $\{\xi_n\}$, $\{\sigma_n\}$ are sequences of positive numbers converging to zero, verifying $\varepsilon_k<1/k^3$, $\sum_{k=1}^\infty \varepsilon_k <\delta$, and
\begin{equation}\label{hombrelobo}
\sqrt{\left(b_k+\frac1{\kappa_2(\partial E^k)}\right)^2+\left(2 \sigma_{k+1}+ b_k\right)^2}-\frac1{\kappa_2(\partial E^k)}+\varepsilon_{k+1}=\frac{\mu'_k}2;
\end{equation}
Furthermore, the sequence $X_n:\overline{\intc P_n} \rightarrow \r^3$ must verify the following properties:
\begin{enumerate}[(A$_{n}$)]
\item $ \sigma_n < \dist_{(\overline{\intc P_n^{\xi_n}},\metri{X_n})}(P_{n-1}^{\xi_{n-1}},P_n^{\xi_n})$;
\item $\| X_n(z)-X_{n-1}(z)\| < \varepsilon_n$, $\forall z\in\intc P_{n-1}^{\varepsilon_{n}}$;
\item $\lambda_{X_n}(z) \geq \alpha_n \lambda_{X_{n-1}}(z)$, $\forall z\in \intc P_{n-1}^{\xi_{n-1}}$, where $\{\alpha_i\}_{i \in \n}$ is a sequence of real numbers such that $0<\alpha_i<1$ and $\{ \prod ^n_{i=1} \alpha_i\}_n$ converges to $1/2$;
\item $\overline{\intc P_{n-1}^{\xi_{n-1}}}\subset \intc P_{n-1}^{\varepsilon_n}\subset \overline{\intc P_{n-1}^{\varepsilon_n}}\subset \intc P_n^{\xi_n} \subset \overline{ \intc P_n^{\xi_n}} \subset \intc P_n\subset  \overline{\intc P_n}\subset \intc P_{n-1}$;
\item $X_n(z)\in E^{n}\setminus (E^{n})_{-b_n}$, for all $z\in P_n$;
\item $X_n(z)\in \r^3\setminus (E^{n-1})_{-2(b_{n-1}+b_n)}$, for all $z\in \intc P_n\setminus \intc P_{n-1}^{\varepsilon_n}$.
\end{enumerate}

The sequence $\{\chi_n\}$ is constructed in a recursive way. To define $\chi_1$, we take $X_1=\varphi$ and $P_1=\Pi$. 

Suppose that we have $\chi_1, \ldots, \chi_n$. In order to construct $\chi_{n+1}$, we consider the following data:
$$E=E^n,\quad E'=E^{n+1},\quad a=b_{n}, \quad X=X_n,\quad P=P_n.$$
Furthermore, Property (E$_n$) says us that $X(P)\subset E\setminus E_{-a}$. Then it is straightforward that we can find a small enough positive constant $\varkappa$, such that Lemma \ref{lem:nadi} can be applied to the aforementioned data, and for any $\varepsilon \in ]0,\varkappa[$.

Take a sequence $\{\widehat{\varepsilon}_m \}\searrow 0$, with $\widehat{\varepsilon}_m<\mbox{minimum}\{\frac{1}{(n+1)^3},\varkappa , b_{n+1}\}$, $\forall m$. For each $m$, we consider $Q_m$ and $Y_m: \overline{\intc Q_m} \rightarrow \r^3$  given by Lemma \ref{lem:nadi}, for the above data and $\varepsilon=b=\widehat{\varepsilon}_m$.
If $m$ is large enough, Assertions (b.1) and (b.5) in Lemma \ref{lem:nadi} tell us that $\overline{\intc P_n^{\xi_{n}}} \subset \intc Q_m$ and the sequence $\{Y_m\}$ converges to $X_n$ uniformly in $\overline{ \intc P_n^{\xi_{n}}}$. In particular, $\{\lambda_{Y_m}\}$ converges uniformly to $\lambda_{X_n}$  in $\overline{ \intc P_n^{\xi_{n}}}$. Therefore there is a $m_0 \in \n$ such that:
\begin{eqnarray} 
\overline{ \intc P_n^{\xi_{n}}} \subset &  \intc P_n^{\hat \varepsilon_{m_0}}& \subset \intc Q_{m_0}, \label{clavao1}\\
\lambda_{Y_{m_0}}& \geq & \alpha_{n+1} \lambda_{X_n} \qquad \hbox{in } \intc P_n^{\xi_{n}}.\label{lambdas}
\end{eqnarray}
We define $X_{n+1}=Y_{m_0}$, $P_{n+1}=Q_{m_0}$, and $\varepsilon_{n+1}=\widehat{\varepsilon}_{m_0}$. The term $\sigma_{n+1}$ is given by (\ref{hombrelobo}) for $k=n+1$. From (\ref{clavao1}) and Statement (b.2), we infer that $\sigma_{n+1} < \dist_{\left(\overline{ \intc P_{n+1}}, \metri{X_{n+1}} \right)}(P_n^{\xi_n},P_{n+1})$. Finally, take $\xi_{n+1}$ small enough such that (A$_{n+1}$) and (D$_{n+1}$) hold. 
The remaining properties directly follow from (\ref{clavao1}), (\ref{lambdas}) and Lemma \ref{lem:nadi}. This concludes the construction of the  sequence $\{\chi_n\}_{n \in \n}$.
\bigskip

Now, we extract some information from the properties of $\{\chi_n\}$. First, from (B$_n$), we deduce that $\{X_n\}$ is a Cauchy sequence, uniformly on compact sets of $\Omega=\bigcup_n \intc P_n^{\varepsilon_{n+1}}=\bigcup_n \intc P_n^{\xi_{n}}$, and so $\{X_n\}$ converges on $\Omega$. If one employs the properties (D$_n$), then the set $\Omega$ is an expansive union of simply connected domains, resulting in $\Omega$ being simply connected. Let $\widetilde \varphi:\Omega\rightarrow \r^3$ be the limit of $\{X_n\}$. Then $\widetilde\varphi$ has the following properties:
\begin{itemize}
\item $\widetilde\varphi$ is a conformal minimal immersion, (Properties  (C$_n$));
\item $\widetilde\varphi:\Omega \longrightarrow D'$ is proper. Indeed, consider a compact subset $K \subset D'$. Let $n_0$ be a natural so that $$K \subset (E^{n-1})_{-2(b_{n-1}+b_n)-\sum_{k\geq n} \varepsilon_k}, \quad \forall n  \geq n_0.$$ From Properties (F$_n$), we have $X_n(z)\in \r^3\setminus (E^{n-1})_{-2(b_{n-1}+ b_n)}$, $\forall z \in \intc P_n \setminus \intc P_{n-1}^{\varepsilon_n}$. Moreover, taking into account (B$_k$), for $k \geq n$, we obtain 
$$\widetilde \varphi(z)\in \r^3\setminus (E^{n-1})_{-2(b_{n-1}+ b_n)-\sum_{k\geq n} \varepsilon_k}.$$
Then, we have $\widetilde\varphi^{-1}(K)\cap(\intc P_n \setminus \intc P_{n-1}^{\varepsilon_n})=\emptyset$ for $n\geq n_0$. This implies that $\widetilde\varphi^{-1}(K)\subset \intc P_{n_0-1}^{\varepsilon_{n_0}}$, and so it is compact in $\Omega$.
\item $\Omega$ is complete with the metric $\metri{ \widetilde{\varphi}}$.  From (\ref{hombrelobo}), $\sigma_n$ is expressible as:
\begin{multline}\label{lobos}
2 \sigma_n=-b_{n-1}+\left(\left(\frac{\mu'_{n-1}}2\right)^2+(\varepsilon_{n})^2+\right.
\\
\left. \frac{\mu'_{n-1}}{\kappa_2(\partial (E^{n-1}))}-\frac{2\varepsilon_{n}}{\kappa_2(\partial (E^{n-1}))}-\frac{4\,b_{n-1}}{\kappa_2(\partial (E^{n-1})))}-\mu'_{n-1}\varepsilon_{n}-4\,b^2_{n-1}\right)^{1/2}.
\end{multline}
Note that  $\{\kappa_2(\partial (E^n))\}$ is a bounded sequence. Then  (\ref{lobos}) lead us to:
$$2 \sigma_n>-\tfrac1{(n-1)^3}+\left(\tfrac{\cte}{{(n-1)}^2}-\tfrac{\cte}{n^3}\right)^{1/2}
>\tfrac{\cte}{n-1} \quad \mbox{($n$ large enough)}.$$
Completeness of $\widetilde{\varphi}$ follows from Properties (A$_n$), (C$_n$), and the fact that $\sum_{n\geq 2} \sigma_n=\infty.$
\item Statement {(a)} in the theorem is a direct consequence of Properties (B$_n$) and the fact $\sum_{n=1}^\infty \varepsilon_n <\delta$.
\item In order to prove Statement {(b)}, we consider $z \in \Omega\setminus \intc \Pi$. There exists $n \in \n$ such that $z \in \intc P_n \setminus \intc P_{n-1}^{\varepsilon_n}.$ Hence, using Properties (B$_k$), $k \geq n$, and Property (F$_n$) we deduce $$\widetilde \varphi(z) \in \r^3 \setminus E^{n-1}_{-2(b_{n-1}+ b_n)-\sum_{k\geq n} \varepsilon_k}\subset \r^3 \setminus D_{-2 d-\delta}.$$
\end{itemize}
 \end{proof}
 Notice that the regularity of $D'$ is used in a strong way to prove the above theorem. Especially, we have made use of this fact to get the lower bound in (\ref{lobos}). Curiously, Theorem \ref{th:mari} can be applied to eliminate the regularity from the statement. So, we are able to establish the following general version.
\begin{teorema} \label{th:refinitivo}
Let $D'$ be an arbitrary convex domain of space. Consider $\psi:\d \longrightarrow \r^3$ a conformal minimal immersion and $\Gamma \subset \d$ a Jordan curve such that $\psi(\Gamma) \subset \partial D$, where $D\subset D'$ is a bounded convex regular domain. Then, for any $\delta>0$, there exist a simply connected domain $\Omega$, with $\Gamma \subset \Omega \subset \d$, and a complete proper minimal immersion $\widetilde{\psi}:\Omega \longrightarrow D'$ so that
\begin{enumerate}[\rm (a)]
\item $\| \widetilde \psi - \psi \| < \delta$, in $\intc \Gamma$;
\item $\widetilde \psi (\Omega \setminus \intc \Gamma) \subset D' \setminus D_{-2 \delta}.$  
\end{enumerate}

 \end{teorema}

\begin{proof}
Without loss of generality, we can assume that $0 \in D'$. If $B(p,r)$ denote the open ball centered at $p$ of radius $r$, then we label $V_n=B(0,n) \cap \left((1-1/n)\cdot D'\right)$, for $n \geq 2$. Notice that $V_n$ is a bounded convex domain. By using Minkowski's theorem (see Theorem \ref{th:minko} in Section \ref{sec:pre}), one has the existence of a decreasing sequence of convex bodies $\{V_n^k \}_k\searrow V_n$,  such that $\partial V_n^k$ is an  analytic surface. As $\dist(\partial D', V_n)>\dist(\partial D', V_{n+1})>0$, we can guarantee the existence of $k_n \in \n$ satisfying $V_n \subset V_n^{k_n} \subset V_{n+1}$. Then we have constructed a increasing sequence of convex domains $D^n \df V_n^{k_n}$ with analytic boundary, verifying:
\begin{itemize}
\item $\cup_{n=2}^{\infty} D^n=D'$;
\item $\overline{D^n} \subset D^{n+1}.$
\end{itemize}
It is clear that we can take a subsequence (denoted in the same way) such that $\overline{D} \subset  D^2$. Choose $d_1>0$ satisfying:  $\overline{D_{d_1/2}} \subset D^2$, $d_1 < \delta/2$, and 
$$\sqrt{\left( d_1 + \frac1{\kappa_2(\partial D_{d_1/2})}\right)^2+{d_1}^2}-\frac1{\kappa_2(\partial D_{d_1/2})}< \frac{\dist (\partial D_{d_1/2}, \partial D^2)}2,$$
and define $D^1 \df D_{d_1/2}.$
At this point, we can also construct a sequence of positive real numbers $\{ d_n\}\searrow 0$ satisfying:
$$\sqrt{\left( d_n + \frac1{\kappa_2(\partial D^n)}\right)^2+{d_n}^2}-\frac1{\kappa_2(\partial D^n)}< \frac{\dist (\partial D^n, \partial D^{n+1})}2.$$
Consider a polygon $\Pi_1$ so that $\psi(\Pi_1) \subset D^1\setminus (D^1)_{-d_1}$, and $\Gamma\subset \intc \Pi_1$. Label $\varphi_1=\psi_{|\overline{\intc \Pi_1}}$. By successive application of Theorem \ref{th:mari}, we will construct a sequence of polygons $\{\Pi_n \}$, a sequence of positive reals $\{ \delta_ n\}\searrow 0$, and a sequence of conformal minimal immersions $\varphi_n: \overline{\intc \Pi_n} \rightarrow D^n$, verifying:
\begin{enumerate}[($\cal A_{\rm n}$)]
\item $ \dist_{(\intc \Pi_n,\varphi_n)}(\Pi_{n},0)>n$;
\item $ \| \varphi_n(z)-\varphi_{n-1}(z) \| < \delta_n ,$ $\forall z \in  \intc \Pi_{n-1},$ where $\displaystyle \sum_{n=1}^\infty \delta_n<\delta$;
\item $ \lambda_{\varphi_n}(z)\geq \alpha_n \lambda_{\varphi_{n-1}}(z),$ $\forall z \in  \intc \Pi_{n-1}$, where  $0<\alpha_i<1$ and $\{ \prod ^n_{i=1} \alpha_i\}_n$ converges to $1/2$;
\item $ \overline{\intc \Pi_{n-1}} \subset \intc \Pi_n$;
\item $ \varphi_n(\intc \Pi_n) \subset D^n \setminus (D^n)_{-d_n}$;
\item $ \varphi_n(\intc \Pi_n \setminus \intc \Pi_{n-1}) \subset \r^3 \setminus (D^{n-1})_{-2 d_n-\delta_n}.$
\end{enumerate}
Assume we have constructed $\varphi_n$. To obtain $\varphi_{n+1}$, we apply Theorem \ref{th:mari} to the following data:
$$ \widehat D=D^n, \quad \widehat D'=D^{n+1}, \quad \widehat \Pi=\Pi_n, \quad \widehat \varphi=\varphi_n, \quad \widehat d=d_n, $$
and $\widehat \delta<\delta_n$ small enough so that ($\cal C_{\rm n+1}$) holds. Then, $\varphi_{n+1}$ is the immersion $\widetilde \varphi$ given by Theorem \ref{th:mari}, and $\Pi_{n+1}$ is a suitable polygon satisfying ($\cal A_{\rm n+1}$), ($\cal D_{\rm n+1}$), and ($\cal E_{\rm n+1}$). Notice that ($\cal B_{\rm n+1}$) and ($\cal F_{\rm n+1}$) directly follows from  Theorem \ref{th:mari}.

From Properties ($\cal B_{\rm n}$), we know that the sequence $\{\varphi_n \}$  converges uniformly on compact sets of $\displaystyle \Omega= \cup_{n=1}^\infty \intc \Pi_n$ (note that $\Omega$ is a simply connected domain in $\d$.) Thus, the uniform limit $\widetilde \psi: \Omega \rightarrow B$ is the immersion we are looking for.

\end{proof}
A trivial consequence of the above theorem is the following corollary.
\begin{corolario} \label{co:esperanza}
For any  convex domain $B$ in $\r^3$ (including noncompact and nonregular ones) there is a proper complete  minimal immersion $\psi: \d \longrightarrow B$.
\end{corolario}
On the other hand, as we mentioned in the introduction, these results give us some extra information about universal regions for minimal surfaces. 
\begin{corolario} \label{co:caridad}
A convex domain of $\r^3$  is not  a universal region for minimal surfaces.
\end{corolario}

 \vspace{10mm}

\noindent Francisco Martín, Santiago Morales \\ Departamento de Geometría y Topología \\  Universidad de Granada \\ 18071 Granada, Spain  \\
\texttt{fmartin@ugr.es, santimo@ugr.es}
\end{document}